\documentclass[11pt]{amsart}

\usepackage{texdraw,multicol,rotating}
\usepackage{epsfig}
\usepackage{amssymb,amsthm,amsfonts,amscd}
\usepackage{latexsym}
\usepackage{fullpage}
\usepackage{amsmath}
\usepackage[vcentermath,enableskew]{youngtab}
\usepackage{graphicx}
\usepackage{stmaryrd}
\usepackage[latin1]{inputenc}
\usepackage[bookmarksnumbered]{hyperref}
\usepackage[arrow, matrix, curve]{xy}

\def\Z{\mathbb{Z}}
\def\R{\mathbb{R}}

\def\Q{\mathbb{Q}}
\def\N{\mathbb{N}}

\def\f{\tilde{f}}
\def\e{\tilde{e}}
\def\ep{\epsilon}
\def\ph{\varphi}
\def\B{\textbf{B}}
\def\rr{\right\rangle}
\def\ll{\left\langle}
\def\#{\sharp}

\def\L{\Lambda}

\def\bea{\begin{eqnarray*}}
\def\eea{\end{eqnarray*}}

\def\M{\mathfrak{M}}
\def\N{\mathfrak{N}}

\def\<#1,#2>{\langle\,#1,\,#2\,\rangle}

\newtheorem{Example}{Example}[section]
\newtheorem{Prop}{Proposition}[section]
\newtheorem{Def}{Definition}[section]

\newtheorem{Lemma}{Lemma}[section]
\newtheorem{Remark}{Remark}[section]
\newtheorem{Cor}{Corollary} [section]

\newtheorem{Theo}{Theorem} [section]

\def\qed{\quad\hfill\ensuremath{\Box}}

\title[Compression of Nakajima monomials \\ in type A and C] {Compression of Nakajima monomials \\ in type A and C}
\author{Matthias Meng$^{*}$}
\date{}
\thanks{$^{*}$Supported by the DFG Graduiertenkolleg 1269  "Globale Strukturen in Geometrie und Analysis" at the University of Cologne}

\begin{document}
\maketitle 

\begin{abstract}
We describe an explicit crystal morphism between Nakajima monomials and monomials which give a realization of crystal bases for finite dimensional irreducible modules over the quantized enveloping algebra for Lie algebras of type A and C. This morphism provides a connection between arbitrary Nakajima monomials and Kashiwara-Nakashima tableaux. This yields a translation of Nakajima monomials to the Littelmann path model. Furthermore, as an application of our results we define an insertion scheme for Nakajima monomials compatible to the insertion scheme for tableaux.
\end{abstract}


\section*{Introduction}

Crystal basis theory for integrable modules over quantum groups as introduced by Kashiwara $[4]$ leads to a combinatorial interpretation of those modules in terms of crystals themselves, and furthermore their various models. Let us list some of those models which will play a role in the present article:
	
\begin{enumerate}
	\item semistandard Young tableaux and reversed Young tableaux, satisfying certain conditions, for classical Lie algebras by Kashiwara and Nakashima $[6]$, and Kim and Shin $[7]$ (see also Kang, Kim, Lee, and Shin in $[2]$ and $[3]$), respectively,
	\item Young walls for affine Lie algebras by Kang, Kim, and Lee $[1]$,
	\item monomials for Kac-Moody algebras discovered by Nakajima  $[11]$, and generalized by Kashiwara  $[5]$,
	\item the path model for symmetrizable Kac-Moody algebras introduced by Littelmann $[10]$.
\end{enumerate}

Let us be more precise about the monomial and the path model. Defining a $t$-analog of $q$-characters Nakajima $[11]$ introduced a set of monomials $\M$ in certain variables $Y_i(k)$, and discovered a crystal structure on certain subsets of $\M$. Kashiwara $[5]$ generalized this, in that he defined a crystal structure on $\M$, and proved that the connected component containing a highest weight monomial of integral weight $\lambda$ is isomorphic to the crystal basis $B(\lambda)$ of irreducible highest weight modules.

Kang, Kim, Lee, and Shin $[2],[3]$ considered specific highest weight monomials $M_\lambda \in \M$ of weight $\lambda$ and gave an explicit description of their connected components $\mathcal{M}(\lambda)$. Furthermore, they exhibited a connection between those and reversed Young tableaux.

As a generalization of Young tableaux Littelmann $[9]$ considered paths (modulo reparametrization) on the real form of the weight lattice and defined the so-called root operators acting on those paths. With these operators the set of paths $\Pi$ becomes a crystal, and every Young tableau can easily be considered as such a path $[10]$.\\

In this article we describe a translation between the monomial and the path model. That is, we map an arbitrary monomial, not necessarily contained in some $\mathcal{M}(\lambda)$, to a path in $\Pi$ such that our mapping yields a crystal morphism. For example, if the underlying Lie algebra is of type $A_1$ a possible definition of such a map is quite obvious: each monomial $M \in \M$ is of the form $M=Y_1(i_1)^{y(i_1)} \cdots Y_1(i_k)^{y(i_k)}$ where $k \in \mathbb{N}$, $i_1, \ldots , i_k \in\mathbb{Z}$ and $i_1<\ldots <i_k$, and $y(i_j)\in \mathbb{Z}$. To a fixed monomial $M$ we associate the path $\pi_M=\pi_{y(i_k)\Lambda_1} * \ldots * \pi_{y(i_1)\Lambda_1}, $ where $\pi_\lambda(t)=t \lambda$ is the path connecting the origin to $\lambda$.

\bigskip

\textbf{Example.} 
\textit{Consider the monomial $M=Y_1(2)^{-1}Y_1(1)^2$:}
\[ 
\xy 
(40,70)*{\scriptstyle\pi};
(50,70)*{}; (90,70)*{} **\dir{.};
(50,73)*{\scriptstyle -\alpha_1};
(60,73)*{\scriptstyle -\Lambda_1};
(70,73)*{\scriptstyle 0};
(80,73)*{\scriptstyle \Lambda_1};
(90,73)*{\scriptstyle \alpha_1};
(50,70)*{\scriptstyle |};
(60,70)*{\scriptstyle |};
(70,70)*{\scriptstyle |};
(80,70)*{\scriptstyle |};
(90,70)*{\scriptstyle |};
(0,70)*{\scriptstyle M=Y_1(2)^{-1}Y_1(1)^2};
(25,70)*{\scriptstyle\longleftrightarrow};
 {\ar (70,70)*{}; (60,70)*{}}; 
 {\ar@/^0.3pc/(60,70)*{}; (80,70)*{}};  
(0,63)*{\scriptstyle\downarrow};
(3,63)*{\scriptstyle\f_1};
(70,63)*{\scriptstyle\downarrow};
(73,63)*{\scriptstyle\f_1};
(40,56)*{\scriptstyle\f_1(\pi)};
(50,56)*{}; (90,56)*{} **\dir{.};
(50,59)*{\scriptstyle -\alpha_1};
(60,59)*{\scriptstyle -\Lambda_1};
(70,59)*{\scriptstyle 0};
(80,59)*{\scriptstyle \Lambda_1};
(90,59)*{\scriptstyle \alpha_1};
(50,56)*{\scriptstyle |};
(60,56)*{\scriptstyle |};
(70,56)*{\scriptstyle |};
(80,56)*{\scriptstyle |};
(90,56)*{\scriptstyle |};
(0,56)*{\scriptstyle \scriptstyle\f_1 M=Y_1(2)^{-2}Y_1(1)^1};
(25,56)*{\scriptstyle\longleftrightarrow};
 {\ar (70,56)*{}; (50,56)*{}}; 
 {\ar@/^0.3pc/(50,56)*{}; (60,56)*{}};  
\endxy
\] 
Note that, even for type $A_2$, to find such a mapping is by far less obvious.

\bigskip

\textbf{Example.}
\textit{For $\frak{g}$ of type $A_2$ we have $\f_1(Y_1(2)^{-1}Y_1(1)^2)=Y_1(2)^{-2}Y_1(1)Y_2(1)$. Adopting the (obvious) construction in type $A_1$, we would associate the path $\pi$ displayed on the left:}

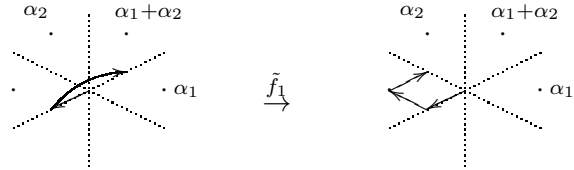
\[ 
\xy 
(20,80)*{}; (20,60)*{} **\dir{.};
(10,75)*{}; (30,65)*{} **\dir{.};
(10,65)*{}; (30,75)*{} **\dir{.}; 
 (30,70)*{\cdot };
 (33,70)*{\scriptstyle\alpha_1 };
  (25,77.5)*{\cdot };
 (28,80.5)*{\scriptstyle\alpha_1+\alpha_2 };
 (15,77.5)*{\cdot };
 (10,70)*{\cdot };
 (13,80.5)*{\scriptstyle\alpha_2};
 {\ar (20,70)*{}; (15,67.5)*{}};
{\ar @/^0.3pc/(15,67.5)*{}; (25,72.5)*{}};
(45,70)*{\scriptsize \stackrel{\f_1}{\rightarrow}};
(70,80)*{}; (70,60)*{} **\dir{.};
(60,75)*{}; (80,65)*{} **\dir{.};
(60,65)*{}; (80,75)*{} **\dir{.}; 
 (80,70)*{\cdot };
 (83,70)*{\scriptstyle\alpha_1 };
  (75,77.5)*{\cdot };
 (78,80.5)*{\scriptstyle\alpha_1+\alpha_2 };
 (65,77.5)*{\cdot };
 (60,70)*{\cdot };
 (63,80.5)*{\scriptstyle\alpha_2};
 {\ar (70,70)*{}; (65,67.5)*{}};
{\ar (65,67.5)*{}; (60,70)*{}};
{\ar (60,70)*{}; (65,72.5)*{}};

\endxy
\] 
\textit{After applying $\f_1$, observe that $\f_1(\pi)$ has a linear part different from any fundamental root direction. That is, $\f_1(\pi)$ does not coincide with the path we would associate to the monomial $\f_1(M)$ in the same manner.}

\bigskip

By generalizing the results of $[2]$ and $[3]$ in type $A$ and $C$ to arbitrary monomials in $\M$, we determine the structure of the crystal graph associated to the connected component of an arbitrary, not necessarily highest weight, monomial in $\M$. More precise, we give a crystal morphism between the set $\M$ and the set of tableaux which give realizations of $B(\lambda)$, and consequently, due to Littelmann $[10]$, we can associate a path to those tableaux.
Our crystal morphism \emph{compresses} an arbitrary Nakajima monomial $M \in \M$ into one which lies in a connected component $\mathcal{M}(\lambda)$, with integral dominant weight $\lambda$ depending on $M$.

In a first step we describe a crystal isomorphism between the Nakajima monomials and certain matrices, namely $\text{Mat}_{n+1 \times \mathbb{Z}} (\Z_{\geq 0})$ in the $A_n$-case and  \\ $\mbox{Mat}_{2n \times \mathbb{Z}} (\Z_{\geq 0})$ in type $C_n$. This bijection allows us to define the compression of 
a monomial by compressing its associated matrix as follows: For simplicity let $M$ denote the matrix associated to an arbitrary monomial $M \in \M$ lying in some a priori unknown connected component of $\M$. We give an algorithm which decomposes $M$ into a sum $M=M_1+M_2$, such that $M_1$ corresponds to a monomial in some $\mathcal{M}(\mu_1)$. Then, we move every column of $M_2$ one step to the left and denote by $M^{(1)}$ the sum of $M_1$ and the altered counterpart of $M_2$. Our procedure allows an iteration yielding a sequence of matrices $M^{(i)}$. Since $M$ has just finitely many nonzero columns, it is guaranteed that after a finite number of steps our iteration becomes stationary and we obtain a matrix $M^{(k)}$ corresponding to a monomial that lies in some $\mathcal{M}(\mu_k)$. We call $M^{(k)}$ the \emph{compressed version} of $M$. Our algorithm respects the crystal structure, that is we prove:

\bigskip

\textbf{Main Theorem.} \textit{Let $\mathfrak{g}$ be of type $A$ or $C$, and let $M \in \M$ be a Nakajima monomial. Denote by $M^{(k)}$ its compressed version. Then, the map 
$$
\begin{array}{ccccl}
  \kappa: &\M&    \rightarrow  &\bigcup\limits_{\lambda} \mathcal{M}(\lambda)& \\
   &M&       \mapsto    & M^{(k)}& \\
    \end{array}              
$$
is a morphism of crystals. In particular, the connected component of $M$ is isomorphic to the connected component of $\kappa(M)$.} \\

Due to $[2]$ and $[3]$ we can assign a tableau $S(N)$ to each $N \in \mathcal{M}(\lambda)$. Consequently, our Main Theorem gives:

\bigskip

\textbf{Corollary.} \textit{Let $\mathfrak{g}$ be of type $A$ or $C$, and let $M \in \M$ be a Nakajima monomial. The mapping sending $M$ to the tableau  $S(\kappa(M))$
yields a crystal morphism.}

\bigskip

Note that, in view of $[10]$ we obtain a translation of Nakajima monomials into Littelmann paths.

\bigskip

\textbf{Example.}
\textit{Consider the monomial $M = Y_1(2)Y_1(1)^2$, and the path obtained via our construction. Observe that our assignment commutes with the crystal operator $\f_1$, as illustrated in the following pictures.}
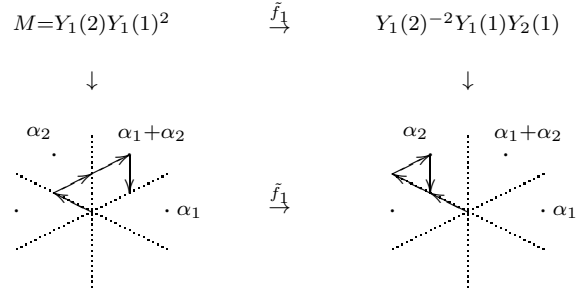
\[ 
\xy 
(20,95)*{\scriptstyle M=Y_1(2)Y_1(1)^2};
(45,96)*{\scriptstyle \stackrel{\f_1}{\rightarrow}};
(70,95)*{\scriptstyle Y_1(2)^{-2}Y_1(1)Y_2(1) };
(20,87.5)*{\scriptstyle \downarrow};
(70,87.5)*{\scriptstyle \downarrow};
(20,80)*{}; (20,60)*{} **\dir{.};
(10,75)*{}; (30,65)*{} **\dir{.};
(10,65)*{}; (30,75)*{} **\dir{.}; 
 (30,70)*{\cdot };
 (33,70)*{\scriptstyle\alpha_1 };
  (25,77.5)*{\cdot };
 (28,80.5)*{\scriptstyle\alpha_1+\alpha_2 };
 (15,77.5)*{\cdot };
 (10,70)*{\cdot };
 (13,80.5)*{\scriptstyle\alpha_2};
 {\ar (20,70)*{}; (15,72.5)*{}};
 {\ar (15,72.5)*{}; (20,75)*{}};
{\ar (20,75)*{}; (25,77.5)*{}};
{\ar (25,77.5)*{}; (25,72.5)*{}};
(45,72)*{\scriptstyle \stackrel{\f_1}{\rightarrow}};
(70,80)*{}; (70,60)*{} **\dir{.};
(60,75)*{}; (80,65)*{} **\dir{.};
(60,65)*{}; (80,75)*{} **\dir{.}; 
 (80,70)*{\cdot };
 (83,70)*{\scriptstyle\alpha_1 };
  (75,77.5)*{\cdot };
 (78,80.5)*{\scriptstyle\alpha_1+\alpha_2 };
 (65,77.5)*{\cdot };
 (60,70)*{\cdot };
 (63,80.5)*{\scriptstyle\alpha_2};
 {\ar (70,70)*{}; (65,72.5)*{}};
{\ar (65,72.5)*{}; (60,75)*{}};
{\ar (60,75)*{}; (65,77.5)*{}};
{\ar (65,77.5)*{}; (65,72.5)*{}};

\endxy
\] 

\bigskip

As another application of our compression and the Corollary we define an insertion scheme for Nakajima monomials compatible with the insertion scheme of reversed tableaux described in $[8]$. More precise, let $M_1$ and $M_2$ be two matrices which correspond to arbitrary monomials in $\M$. Then, we consider the matrix $M_1*M_2=(M_2,\mathbf{0},M_1)$ with a suitable zero-matrix $\mathbf{0}$ and apply our compression procedure to $M_1*M_2$. Following the convention that $M_1*M_2$ interchangebly denotes the matrix and its associated monomial, we obtain $\kappa(M_1*M_2)\in \bigcup_{\lambda} M(\lambda)$ and the tensor product rule of crystals yields

\bigskip

\textbf{Theorem.} 
{\textit{Let $\mathfrak{g}$ be of type $A$ or $C$, and let $\M$ be the set of Nakajima monomials. Then, the map 
 $$
\begin{array}{ccccl}
  &\M \otimes \M&    \rightarrow  &\bigcup\limits_{\lambda} \mathcal{M}(\lambda)& \\
   &M_1\otimes M_2&       \mapsto    &\kappa(M_1*M_2)& \\
    \end{array}              
$$
is a morphism of crystals. In particular, the connected component of $M_1\otimes M_2$ is isomorphic to the connected component of $\kappa(M_1*M_2)$.}

\section{Nakajima monomials}

In this section we define the Nakajima monomials and their crystal structure. 
Let $\frak{g}$ be  an arbitrary symmetrizable Kac-Moody Lie algebra with weight lattice $P$ and $I$ an index set such that $\alpha_i \in P$ for $i \in I $ are the simple roots. Let further $h_i \in P^*$ be the simple coroots and $(\cdot , \cdot):P \times P \rightarrow \Q$ a bilinear symmetric form. For $ i \in I $ and $\lambda \in P$ set $\left\langle h_i,\lambda \right\rangle  := \frac{2(\alpha_i,\lambda)}{(\alpha_i,\alpha_i)}$.  \

For $i\in I$ and $n\in \Z$ we consider monomials in the variables $Y_i(n)$. That means we obtain the set of Nakajima monomials $\M$ as follows

\begin{center}

 $ \M:=\bigg\{\prod\limits_{i \in I,n\in \Z} Y_i(n)^{y_i(n)}  ;  y_i(n) \in\Z $ vanish except for finitely many $(i,n)\bigg\}.$ 
 \end{center} 

In order to define the crystal structure on $\M$ we take some integers \\ $c= (c_{ij})_{i \neq j \in I} \subset \Z$  such that $c_{ij}+c_{ji}=1$ and consider the monomials 
\begin{center}

$ A_i(n):= Y_i(n)Y_i(n+1) \prod\limits_{j \neq i} Y_j(n+c_{ji})^{\left\langle h_j,\alpha_i \right\rangle} .$
                                   
\end{center}

Let now $M$ be an arbitrary monomial in $\M$ and $i\in I$. Then we set: 

\begin{center} 
  $
    \begin{array}{ccl}    wt(M)    &=&     \sum\limits_{i}(\sum\limits_{n} y_i(n))\L_i, \\
                         \ph_i(M)  &=&         \mbox{max} \{ \sum\limits_{k\leq n} y_i(k) ; n \in \Z\}, \\             
                         \ep_i(M)  &=&       \mbox{max}   \{-\sum\limits_{k > n} y_i(k) ; n \in \Z\}, \\
    \end{array}
  $
    
\end{center}
where $\L_i \in P $ are the fundamental weights, that means $\left\langle h_j,\L_i \right\rangle= \delta_{i,j}$.
To define the operators $\e_i$ and $\f_i$ we consider the values

\begin{center}
       $
     \begin{array}{ccl}  n_f  &=&  \mbox{min} \{n;\ph_i(M)=\sum\limits_{k\leq n} y_i(k)\} \\  
                              &=&    \mbox{min} \{n; \ep_i(M)= -\sum\limits_{k > n} y_i(k)\},\\ \\
                             
                         n_e  &=& \mbox{max} \{n;\ph_i(M)=\sum\limits_{k\leq n} y_i(k)\} \\  
                              &=&    \mbox{max} \{n; \ep_i(M)= -\sum\limits_{k > n} y_i(k)\}\\ \\
     \end{array}
     $
\end{center} 
and set
\begin{center}
      $
      \begin{array}{ccl}
         \f_i(M)  &=& \begin{cases} 
                                                       0                & \text{if   $ \ph_i(M)=0,$} \\
                                               A_i(n_f)^{-1} M   &\text{if    $\ph_i(M) > 0,$}
                              \end{cases}                 
                                                \\ \\

                            \e_i(M)  &=& \begin{cases} 
                                                       0          & \text{if  $  \ep_i(M)=0,$} \\
                                               A_i(n_e) M   & \text{if  $    \ep_i(M)>0.$} 
                                              \end{cases}
      \end{array}   
      $
\end{center} 

\begin{Prop}$[5]$
With the maps $wt, \ph_i,\ep_i,\f_i$ and $\e_i$ thus defined, $\M$ becomes a semi-normal crystal.

\end{Prop}

We denote this crystal by $\M_c$ because the crystal structure of $\M$ depends on the choice of $c$. On the other hand one can easily see that the isomorphism class of the crystal $\M_c$ does not depend on this choice.  \\
From now on, for simplicity, we choose $ c= (c_{ij})_{i \neq j \in I}$ as follows: 

\begin{center} 
$c_{ij}= \begin{cases}     0    &\text{if  $i>j$,} \\      
                                       1    &\text{else.}  
              \end{cases}                         
$               
\end{center}

Now we recall the following result of Kashiwara.

\begin{Prop}  $[5]$
Let $M$ be a monomial of weight $\lambda$ with $\e_i(M)=0$ for all $i\in I$. Then the connected component of $\M$ containing $M$ is isomorphic to $B(\lambda)$.  
\end{Prop}

The aim of this thesis is to give such an isomorphism explicitly for not necessarily highest weight monomials. 
In the first part we define this isomorphism for Lie algebras of type $A$. In the second part we generalize this to type $C$.

\section{Compression of Nakajima monomials in type $A$}
Henceforth we consider a Lie algebra $\frak{g}$ of type $A_n$. In this case we have the fundamental weights $\L_1, \ldots, \L_n$ and we get an orthogonal basis $\beta_1, \ldots ,\beta_{n+1}$ with $\beta_1=\L_1, \beta_i=\L_i-\L_{i-1}$ for $2 \leq i  \leq n$ and $\beta_{n+1}=-\L_n$. Moreover the simple roots are given by $\alpha_i=\beta_i-\beta_{i+1}$. 
Thus we compute
\begin{center}
$A_i(j)=Y_i(j)Y_i(j+1)Y_{i-1}(j+1)^{-1}Y_{i+1}(j)^{-1}.$
\end{center}
For $i\in \{1,\ldots ,n+1\}$ and $j \in \Z$ we introduce some specific monomials which will be of special interest to us:  
\begin{center}
$X_i(j):= Y_{i-1}(j+1)^{-1}Y_i(j),$
\end{center}
where we set $Y_{n+1}(j)=1=Y_0(j)$ for all $j \in \Z$. \\
With this notation we observe:
\begin{center}
$A_i(j)=X_i(j)X_{i+1}(j)^{-1}.$
\end{center}
Let us briefly recall the monomial realization of the crystal bases $B(\lambda)$ given in $[2]$:

\begin{Prop} $[2]$
Let $\lambda=\sum\limits_{k=1}^{n}a_k\Lambda_k$ be a dominant integral weight and consider $M_1=Y_1(1)^{a_1}Y_2(1)^{a_2} \ldots Y_n(1)^{a_n}$ as highest weight monomial. Then the connected component $\mathcal{M}_1(\lambda)$ of $\M$ containing $M_1$ is characterized as the set of monomials of the form

\begin{center}

$M=\prod\limits_{{i\in \{1, \ldots,n+1 \}} \atop{j \in \{1, \ldots,n \}} } X_i(j)^{m_{ij}}$

\end{center} 
with 

\begin{itemize}
\item[(i)]   $\sum\limits_{i=1}^{n+1}m_{ij}=a_{j+1}+ \ldots + a_n$ for $j=1, \ldots , n,$

\item[(ii)]  $\sum\limits_{k=i}^{n+1}m_{k,j} \leq \sum\limits_{k=i+1}^{n+1}m_{k,j-1}$ for $j=2, \ldots , n+1$ and $i=1, \ldots , n+1.$
\end{itemize}

\end{Prop}
For $s\in \mathbb{Z}$ we also consider the following \emph{shifted} highest weight monomials of weight $\lambda$ $$M_s=Y_1(s)^{a_1}Y_2(s)^{a_2} \ldots Y_n(s)^{a_n}.$$ 
As an immediate consequence of Proposition 2.1 we obtain their connected component $\mathcal{M}_s(\lambda)$ by the set of monomials of the form 
$$M=\prod\limits_{{i\in \{1, \ldots,n+1 \}} \atop{j \in \{s, \ldots,s+n-1\}} } X_i(j)^{m_{ij}}$$
satisfying condition $(i)$ for $j=s, s+1, \ldots,s+n-1$  and $(ii)$ for $i=1, \ldots , n+1$ and $j=s+1,\ldots, s+n.$  \\

Our aim is to \emph{compress} an arbitrary monomial into the form of those in $\mathcal{M}_s(\lambda)$ for a suitable $\lambda \in P$ and $s\in \mathbb{Z}$ such that the crystal structure is preserved. \\ As a first step we write an monomial in $\M$ as a product of $X_i(j)$`s. Thus we show that $\M$ is generated by the variables $X_i(j)$. That means we consider $\M$ as a group with the multiplication of monomials as binary operation. Let $\mathbf{M}$ be the free abelian monoid generated by the set $\{X_i(j),i\in \{1,\ldots,n+1\},j\in \Z\}$, with the same operation and we define an ideal $\mathbf{J}\subset \mathbf{M}$ by
$$
\mathbf{J}=\langle \prod\limits_{k=1}^{n+1} X_k(j+i-k), \mbox{ for } i=1,\ldots,n+1 \mbox{ and } j\in\Z  \rangle_{\mathbf{M}}.
$$ 
The quotient $\mathbf{M}/\mathbf{J}$ becomes a group since we obtain the inverse of $X_i(j)$ by $\prod\limits_{k=i+1}^{n+1} X_k(j-k+i)\prod\limits_{k=1}^{i-1} X_k(j+i-k).$ Moreover we get

\begin{Prop} Sending $X_i(j)$ onto $Y_i(j)Y_{i-1}(j+1)^{-1}$ yields a group isomorphism and therefore we get
$$\M\cong \mathbf{M}/\mathbf{J}.$$
\end{Prop}

\begin{proof} In order to show surjectivity let $M$ be of the form $M=\prod\limits_{i\in I,j \in \Z} Y_i(j)^{y_i(j)}$.  \\
First we write every $Y_{i-1}(j+1)^{-1}Y_i(j)$ that already occurs in $M$ as $X_i(j)$. Then we consider the other $Y_i(j)^{y_i(j)}$`s in $M$. There are two possible cases:  \\
\underline{1. case}: $y_i(j)>0$. Then we write
\begin{center}
$
\begin{array}{ccl}
Y_i(j)&=&\prod\limits_{k=0}^{i-1} Y_k(j+i-k)^{-1} \prod\limits_{k=1}^{i}Y_k(j+i-k)\\ &=& \prod\limits_{k=1}^{i} X_k(j+i-k).
\end{array}
$
\end{center}
Therefore we get
\begin{center} 
$Y_i(j)^{y_i(j)}= \prod\limits_{k=1}^{i} X_k(j+i-k)^{y_i(j)}$. 
\end{center}
\underline{2. case}: $y_i(j)<0$. In this case we get 
\begin{center} 
$
  \begin{array}{ccl}
Y_i(j)^{-1} &=& \prod\limits_{k=i}^{n} Y_k(j-(k-i))^{-1}\prod\limits_{k=i+1}^{n+1}Y_k(j-(k-i)) \\
            &=& \prod\limits_{k=i+1}^{n+1} X_k(j-k+i).
  \end{array}
$
\end{center}
Hence we have 
\begin{center} 
$Y_i(j)^{y_i(j)}= \prod\limits_{k=i+1}^{n+1} X_k(j-k+i)^{-y_i(j)}$. 
\end{center}
These equations imply 
$$Y_i(j)^{-1}Y_{i-1}(j+1)=\prod\limits_{k=i+1}^{n+1} X_k(j-k+i)\prod\limits_{k=1}^{i-1} X_k(j+i-k)$$
and hence with the definition of $\mathbf{J}$ we have injectivity.

\qed
\end{proof}

Let now $M \in \M$ be an arbitrary monomial. Due to Proposition 2.2 we can write $M$ as a product of $X_i(j)$`s. That means we find $m_{ij} \in \Z_{\geq 0}$ such that 
$$M=\prod\limits_{i\in \{1,\ldots ,n+1\},j \in \Z} X_i(j)^{m_{ij}}.$$

Writing $M$ in this way is obviously not unique. But we can fix a \emph{reduced notation} $[m_{i,j}]$ and associate this matrix. Let us define the reduced notation on the level of matrices.  \\

Let $M=(m_{ij})$ be an arbitrary matrix in $\mbox{Mat}_{n+1 \times \Z} (\Z_{\geq 0})$, where $\mbox{\mbox{Mat}}_{n+1 \times \mathbb{Z}} (\Z_{\geq 0})$ is the set of matrices with  infinitely many columns but just  finitely many different from zero.  Then we get the reduced form $[M]$ of $M$ by applying the following rule:
\begin{itemize}

\item[\textbf{(A1)}] 
For every $i\in \{ 1, \ldots,n+1 \}$ we search for $j\in \Z$ such that 
\begin{center}
 $m_{i+s,j-s} \neq 0 $ for all $s=-i+1,-i+2, \ldots ,-1,0,1,\ldots, n-i,$  
\end{center}
 then we decrease these entries by 
 \begin{center}
 min$\{m_{i+s,j-s}; s=-i+1,-i+2, \ldots ,-1,0,1,\ldots, n-i \}$. 
\end{center} 
\end{itemize}

Denote by $[M]$ the matrix obtained from $M$ by applying this rule.  \\ \\
From now on we associate a matrix to a monomial in the following way: We write every $Y_i(j)^{y_i(j)}$ as a product of $X_k(l)$`s as in Proposition 2.1 and get a corresponding matrix $M=m_{ij}$. Then we apply (A1) and obtain $[M]$.     
We define an equivalence relation on  $\mbox{Mat}_{n+1 \times \mathbb{Z}} (\Z_{\geq 0})$ by 
\begin{center}
$M \sim N$ iff $[M]=[N]$  
\end{center}
and consider the quotient  
\begin{center}
$\mbox{\mbox{Mat}}_{n+1 \times \mathbb{Z}} (\Z_{\geq 0})/ \sim$ .
\end{center}
Now it is obvious that two matrices $M \sim N$ correspond to the same monomial and therefore we get a well defined map by sending a monomial to the associated matrix $[M]$ as above. Moreover we obtain

\begin{Prop}
There exists a bijection between $\M$ and $\text{Mat}_{n+1 \times \mathbb{Z}} (\Z_{\geq 0})/ \sim$.   
\end{Prop}

%
%
%

Since we want this bijection to become a crystal morphism we need to define a crystal structure on $\mbox{\mbox{Mat}}_{n+1 \times \mathbb{Z}} (\Z_{\geq 0})/ \sim$ which coincides with the structure on $\M$ under our bijection.  \\  \\
Let $M=(m_{ij})_{ {i=1,\ldots,n+1} \atop{j\in\Z}} \subset \Z_{\geq 0} $ be a $(n+1) \times \Z$-matrix.  \\
Set 
\begin{center}
$
\begin{array}{ccl}
wt(M) &=& \sum\limits_{i=1}^{n+1}(\sum\limits_{j\in\Z} m_{ij})\beta_i, \\ 
\ph_i(M) &=&  $ max $ \{\sum\limits_{j \leq k}m_{ij} - \sum\limits_{j <k} m_{i+1,j} $ ; $ k \in \Z\}, \\      
 \ep_i(M) &=& -$ min $\{ \sum\limits_{j > k}m_{ij}- \sum\limits_{j \geq k} m_{i+1,j} $ ; $ k \in \Z\}. \\
\end{array}
$
\end{center}
If $\ph_i(M)=0$ we set $\f_i(M)=0$. 
Otherwise let $k \in \Z$ be minimal such that $$ \ph_i(M)=\sum\limits_{j \leq k}m_{ij} - \sum\limits_{j <k} m_{i+1,j}. $$  
Note that, this $k$ exists because $M$ has just finitely many columns different from zero.\\
We define $\f_i(M)$ as the matrix we get from $M$ by increasing (resp. decreasing) $m_{i+1,k}$ (resp. $m_{i,k}$) by one. Formally spoken we get $\f_i(M)=\hat{m}_{s,j} $  from $M=m_{s,j}$ by \\
\begin{center}
            $ \hat{m}_{s,j}  = \begin{cases} 
                                                   m_{s,j}            & \text{if  $(s,j) \notin \{(i,k),(i+1,k)\}$} , \\
                                               m_{i,k}-1              & \text{if  $(s,j)=(i,k)$,} \\
                                               m_{i+1,k}+1           & \text{if $(s,j)=(i+1,k)$.} \\
                                            \end{cases}    
            $
\end{center}                                            

Similarly, we can define the operator $\e_i$ : \\
If $\ep_i(M)=0$ we set $\e_i(M)=0$. \\ 
For $\ep_i(M)\neq0$ let $ p \in \Z$ be maximal such that $$\ep_i(M) = -  (\sum\limits_{j > p}m_{ij}- \sum\limits_{j \geq p} m_{i+1,j}). $$ 
Then we obtain $\e_i(M)=\hat{m}_{s,j}$ from $M=m_{s,j}$ by 

\begin{center}
           $  \hat{m}_{s,j}  = \begin{cases}
                                                   m_{s,j}            & \text{if   $(s,j) \notin \{(i,p),(i+1,p)\}$ ,} \\
                                               m_{i,p}+1            & \text{if   $ (s,j)=(i,p)$,} \\
                                               m_{i+1,p}-1         & \text{if   $(s,j)=(i+1,p)$.} 
                                            \end{cases}    
             $                               
\end{center}      
It is easy to see that these maps are well defined and that  $\mbox{Mat}_{n+1 \times \mathbb{Z}} (\Z_{\geq 0})/ \sim$ along with $wt,\ph_i,\ep_i,\f_i$ and $\e_i$ becomes a semi-normal crystal. \\
Now we can prove:
\begin{Prop} The bijection 
\begin{center}
$\begin{array}{ccccl}
 \Psi:  &\M&    &\rightarrow&  \mbox{Mat}_{n+1 \times \mathbb{Z}} (\Z_{\geq 0})/ \sim   \\
&M=\prod\limits_{i\in \{1,\ldots ,n+1\},j \in \Z} X_i(j)^{m_{ij}}&   &\mapsto&    [m_{ij}]
\end{array}$
\end{center}
is a crystal isomorphism. 
\end{Prop}

\begin{proof} We have to verify that for every $M \in \M$ and $i \in I$ the following holds: 
$$
 wt(M)=wt(\Psi(M)) ,  
$$ 
$$
\ph_i(M)= \ph_i(\Psi(M)), 
$$
$$
     \ep_i(M)= \ep_i(\Psi(M)),
$$
$$    
\Psi(\f_i(M))=\f_i(\Psi(M)),
$$
$$
\Psi(\e_i(M))=\e_i(\Psi(M)).
$$
So let $M=\prod\limits_{s \in I,t \in \Z} Y_s(t)^{y_s(t)}\in \M$ be arbitrary and $\Psi(M) \in \mbox{Mat}_{n+1 \times \mathbb{Z}} (\Z_{\geq 0})$ its corresponding reduced matrix. Now we show that
$$wt(M)=wt(\Psi(M)).$$
Assume we write $M$ as a product of $X_s(t)$`s by writing every factor $Y_s(t)^{y_s(t)}$ as in Proposition 2.1 with corresponding matrix $m_{s,t}$. Now it suffices to show that $wt(M)$ coincides with $wt(m_{s,t})$ because it is obvious that $wt$ is invariant under (A1). So we get \\
\begin{center}
$
\begin{array}{ccccl}
&& wt(M)&=&   \sum\limits_{s}(\sum\limits_{t} y_s(t))\L_s \\  
&\mbox{since }\L_s=\beta_1+\ldots+\beta_s,&   &=& \sum\limits_{s}(\sum\limits_{t} y_s(t))(\beta_1+\ldots+\beta_s) \\
&&    &=& \sum\limits_{s}(\sum\limits_{t}(\underbrace{\sum\limits_{t\leq s} y_s(t))}_{=m_{st}})\beta_s\\
&&    &=& \sum\limits_{s}(\sum\limits_{t}m_{st})\beta_s\\
&\ \mbox{since $wt$ is invariant under (A1),}&    &=& wt(\Psi(M)).

\end{array}
$ 
\end{center}
The same computations work for $\ph_i$ and $\ep_i$. \\ 

Now we show that $\Psi$ commutes with $\f_i$: \\
Let $n_f$ be minimal such that $\ph_i(M)=\sum\limits_{t\leq n_f} y_s(t).$ \\ Then we get $\f_i(M)=A_i(n_f)^{-1}M$.
Due to the choice of $n_f$ we know that $n_f$ is minimal such that\\
$$\ph_i(M)=\ph_i(\Psi(M))=\sum\limits_{j \leq n_f}m_{ij} - \sum\limits_{j <n_f} m_{i+1,j}.$$
That means we decrease (resp. increase)  $m_{i,n_f}$ (resp.  $m_{i+1,n_f}$) by one in $\Psi(M)$. \\
But since $A_i(n_f)^{-1}=X_i(n_f)^{-1}X_{i+1}(n_f)$ it follows that $\f_i(\Psi(M))$ is a corresponding matrix of $\f_i(M)$. It remains to show that $\f_i(\Psi(M))=[\f_i(\Psi(M))].$\\
Assume we had to apply (A1) only after having operated with $\f_i$ but not before. That means we get a full diagonal $m_{i+1+s, n_f-s} \neq 0$ for all $s=-i+1,-i+2, \ldots ,-1,0,1,\ldots, n-i$ after having increased $m_{i+1,n_f}$ by one. But due to the choice of $n_f$ we have $m_{i,n_f+1}<m_{i+1,n_f}$ since otherwise $\sum\limits_{j \leq n_f}m_{ij} - \sum\limits_{j <n_f} m_{i+1,j}$ would not be maximal. Therefore increasing $m_{i+1,n_f}$ doesn`t cause any new (A1) application and 
$$   
\Psi(\f_i(M))=\f_i(\Psi(M)).
$$
The same arguments hold for $\e_i$ which finishes our proof.

 \qed
\end{proof}

Now we define the set of matrices such that the corresponding monomials give a realization of the crystal bases $B(\lambda)$.

\begin{Def}   Define $\N \subset \mbox{Mat}_{n+1 \times \mathbb{Z}} (\Z_{\geq 0})$ as the set of matrices whose reduced forms have only zero-entries out of an $(n+1)\times n$-submatrix $M=(m_{ij})_{{i=1,\ldots,n+1}\atop{j=0,\ldots,n-1}}$ with the following properties:
\begin{itemize}
\item[(i)]  $m_{ij} \in \Z_{ \geq 0}$ for $i=1,\ldots,n+1$ and $j=0,\ldots,n-1$.
\item[(ii)] $\sum\limits_{k=i}^{n+1}m_{k,j} \leq  \sum\limits_{k=i+1}^{n+1}m_{k,j-1}$ for $i=1,\ldots,n+1$ and $j=1,\ldots,n-1$, \\
where we set $\sum\limits_{k=i+1}^{n+1}m_{k,j-1}=0$ for $i=n+1.$
\end{itemize}

\end{Def}  
Due to Proposition 2.1 the associated monomials of matrices in $\N$ can be considered as elements in $\mathcal{M}_s(\lambda)$ for a suitable $\lambda \in P$ and $s\in \mathbb{Z}$. Hence instead of a crystal morphism between $\M$ and $\bigcup_{\lambda\in P,s\in \mathbb{Z}} \mathcal{M}_s(\lambda)$ which we originally intended to find, we just need a morphism from $\mbox{Mat}_{n+1 \times \mathbb{Z}} (\Z_{\geq 0})/ \sim$ to $  \N/ \sim$.  \\ The idea is to compress the matrices. More precisely we move entries into the next column to the left such that the crystal structure is preserved. We do this by decomposing our matrix $M=M_1+M_2$ with $M_1\in \N$ according to the following rule:   \\

\textbf{The lower decomposition rule:} \\
Let $M=[M]=m_{ij}\in \mbox{Mat}_{n+1 \times \mathbb{Z}} (\Z_{\geq 0})$ be a reduced version of an arbitrary matrix. Let $k\in \Z$ be minimal and $l \in \Z$ be maximal such that $m_{ij}=0$ for all $j<l,j>k$ and $i\in \{1,\ldots,n+1\}$. That means the finite part of $M$ which is different from zero is an $(n+1) \times (l-k+1)$-matrix over $\Z_{\geq 0}.$ For simplicity we set $p=l-k$ and  renumber the columns by $0,\ldots,p$. We also assume that $p\geq n-1$, otherwise we fill the matrix with zero-entries on the right side.    \\ We search for $M_1 \in \N$ such that 
$$
M=M_1+M_2.
$$
We explain how to compute $M_1=m_{ij}^{(1)}$ out of $M=m_{i,j}$ recursively: \\
For $i=1,\ldots, n+1$ we set $m_{i,0}^{(1)}:=m_{i,0}$. \\ 
Then, for each $j$ from $1$ to $p$ we do the following: \\
For $i=n+1$ to $i=1$ we compare  
$$\sum\limits_ {k\geq i+1}m_{k,j-1}^{(1)} \ \ \mbox{with} \  \ m_{i,j}+\sum\limits_ {k\geq i+1}m_{k,j}^{(1)}$$
and if $\sum\limits_ {k\geq i+1}m_{k,j-1}^{(1)} < m_{i,j}+\sum\limits_ {k\geq i+1}m_{k,j}^{(1)}$ then we set 
$$m_{i,j}^{(1)} :=  \sum\limits_ {k\geq i+1}m_{k,j-1}^{(1)}-\sum\limits_ {k\geq i+1}m_{k,j}^{(1)}.$$
Otherwise, namely if $\sum\limits_ {k\geq i+1}m_{k,j-1}^{(1)} \geq m_{i,j}+\sum\limits_ {k\geq p+1}m_{k,j}^{(1)},$  we set
$$m_{i,j}^{(1)} = m_{i,j}.$$ 


 This way we get $M_1$ and set $$M_2:=M-M_1.$$
 By construction it is obvious that $M_1$ satisfies condition $(i)$ and $(ii)$ of Definition 2.1 but it remains to show that $M_1$ has at most $n$ columns different from zero such that we can guarantee that $M_1\in \N$. For that we show 
 
 \begin{Lemma}
 Let $[M]=(m_{ij})_{{i=1,\ldots, n+1} \atop {j=0,\ldots l-1}} \in \mbox{Mat}_{n+1 \times l} (\Z_{\geq 0})$ be a reduced matrix without zero columns which satisfies condition $(ii)$ of Definition $2.1$.\\ Then we have $$l\leq n.$$   
 \end{Lemma}
\begin{proof}
Since $\sum\limits_{k=i+1}^{n+1}m_{k,0}=0$ for $i=n+1$ condition $(ii)$ yields $m_{n+1,1}=0$. Again by condition $(ii)$ we obtain $0=m_{n+1,1}\geq m_{n,2}+m_{n+1,2}$ and therefore $m_{n,2}=m_{n+1,2}=0$. 
In general condition $(ii)$ provides 
$$\sum\limits_{k=n+2-j}^{n+1}m_{k,j}=0.$$

That means in particular that $m_{i,n}=0$ for all $i=2,\ldots,n+1$ and $m_{i,j}=0$ for all $i\in I$ and $j\geq n+1$. \\ It remains to show that $m_{1,n}=0$.\\ 
Assume $m_{1,n} \neq 0$. Since $\sum\limits_{k=n+2-j}^{n+1}m_{k,j}=0$ this implies $m_{1+t,n-t}\geq m_{1,n}$ for all $t=1,\ldots ,n$. That means we can apply (A1) which is a contradiction to $M$ being in reduced form and hence 
$$m_{1,n}=0.$$
\qed
\end{proof}

%
%
%
%
%
%
%
%
%
%

Now we can define our desired map:

$$
\begin{array}{ccccccl}
&\Phi&   :   &\mbox{Mat}_{n+1 \times \mathbb{Z}} (\Z_{\geq 0})/ \sim&    \rightarrow  &\N / \sim& \\
     &&         &m_{ij}&    \mapsto       &n_{ij},&
     \end{array}              
$$
where the matrix $n_{ij}$ is computed in the following way: \\ \\
Let $M$ be the reduced version of an arbitrary matrix in $\mbox{Mat}_{n+1 \times \mathbb{Z}} (\Z_{\geq 0})$. Then we consider the lower decomposition of $M$:
$$M=M_1+M_2$$
with $M_1\in \N$. \\
Then we move every entry of $M_2$ one column to the left and denote the new $M_2$ by $M_2^{(1)}.$ Now we set $M^{(1)}:=M_1+M_2^{(1)}$ and if $M^{(1)}\in\N$ we are done and set 
$$M^{(1)}=:N=n_{ij}.$$ 

If  $M^{(1)}\notin\N$ we consider the lower decomposition of $[M^{(1)}]$ and do the same again. This iteration yields a sequence of matrices $M^{(i)}$. Since $M$ has just finitely many columns different from zero there exists a $k\in \mathbb{N}$ such that the iteration becomes stationary with $M^{(k)}\in \N $ and we set: 
$$M^{(k)}=:N=n_{ij}.$$
Let us now combine $\Phi$ and $\Psi$ to obtain the \emph{compression map} $\kappa$ from the set of arbitrary monomials into the set of monomials which give a realization of the crystal bases $B(\lambda)$:

$$
\begin{array}{ccccccl}
\kappa:=&\Psi^{-1} \circ \Phi \circ \Psi &   :   &\M&    \rightarrow  &\bigcup\limits_{\lambda\in P,s\in \mathbb{Z}} \mathcal{M}_s(\lambda).& \\
     \end{array}              
$$
Before we prove that this map is a crystal morphism we consider an example:

\begin{Example}As above, let $\mathfrak{g}$ be of type $A_4$ and take the monomial 
$$M=Y_1(4)^{-1}Y_3(1) Y_1(3)^{-1}Y_4(1)^{-1}Y_2(0)^2 Y_3(2)^2.$$  
Due to Proposition $2.2$ we write:
$$
\begin{array}{ccl}
Y_1(4)^{-1}&=&X_2(3)X_3(2)X_4(1)X_5(0),  \\
Y_3(1)&=&X_3(1)X_2(2)X_1(3),  \\
Y_1(3)^{-1}&=&X_2(2)X_3(1)X_4(0)X_5(-1), \\
Y_4(1)^{-1}&=&X_5(0),  \\
Y_2(0)^2&=&X_2(0)^2X_1(1)^2 ,\\

Y_3(2)^2&=&X_3(2)^2X_2(3)^2X_1(4)^2. \\

\end{array}
$$
That means we obtain the associated reduced matrix by 
$$\begin{pmatrix}     0 & 0 & 2 & 0 &  1 & 2 \\ 
                                    0 & 2 & 0 & 2  & 3 & 0\\       
                                    0 & 0 & 2 & 3  & 0 & 0\\
                                    0 & 1 & 1 & 0  & 0 & 0 \\
                                    1 & 2 & 0 & 0  & 0 & 0

       \end{pmatrix}
       \stackrel{(A1)}{=}
       \begin{pmatrix}    0 & 2 & 0 &  0 & 1 \\ 
                                     2 & 0 & 1  & 2 & 0\\       
                                     0 & 1 & 2  & 0 & 0\\
                                     0 & 0 & 0  & 0 & 0 \\
                                     1 & 0 & 0  & 0 & 0
       \end{pmatrix}=[M],
$$
where we always only consider the finite part of the matrix which is different from zero. \\
We obtained the lower decomposition of $M$ by:
$$
 \begin{pmatrix}          0 & 2 & 0 &  0 & 1 \\ 
                                     2 & 0 & 1  & 2 & 0\\       
                                     0 & 1 & 2  & 0 & 0\\
                                     0 & 0 & 0  & 0 & 0 \\
                                     1 & 0 & 0  & 0 & 0       
\end{pmatrix}=
 \begin{pmatrix}          0 & 2 & 0 &  0 & 0 \\ 
                                     2 & 0 & 1  & 0 & 0\\       
                                     0 & 1 & 0  & 0 & 0\\
                                     0 & 0 & 0  & 0 & 0 \\
                                     1 & 0 & 0  & 0 & 0
       \end{pmatrix}+
 \begin{pmatrix}          0 & 0 & 0 &  0 & 1 \\ 
                                     0 & 0 & 0  & 2 & 0\\       
                                     0 & 0 & 2  & 0 & 0\\
                                     0 & 0 & 0  & 0 & 0 \\
                                     0 & 0 & 0  & 0 & 0       \end{pmatrix}$$ and therefore  
$$M^{(1)}= \begin{pmatrix}      
                                     0 & 2 & 0 &  1  \\ 
                                     2 & 0 & 3  & 0 \\       
                                     0 & 3 & 0  & 0 \\
                                     0 & 0 & 0  & 0  \\
                                     1 & 0 & 0  & 0  
       \end{pmatrix}
$$ 
Since $M^{(1)}\notin \N$ we decompose $[M^{(1)}]=M^{(1)}$ and get 
$$M^{(2)}=
\begin{pmatrix}      
                                     0 & 2 & 0 &  1  \\ 
                                     2 & 2 & 1  & 0 \\       
                                     2 & 1 & 0  & 0 \\
                                     0 & 0 & 0  & 0  \\
                                     1 & 0 & 0  & 0  
       \end{pmatrix} \in \N
$$
Now we apply $\Psi^{-1}$ and receive the monomial          
$$
\begin{array}{ccl}
N=\Psi^{(-1)}(M^{(2)}) &=&  X_1(1)^2X_1(3)X_2(0)^2X_2(1)X_3(0)^2X_3(1)X_5(0) \\
                                      &=& Y_1(3)Y_2(0)^2Y_1(2)^{-2}Y_3(0)^2Y_3(1)Y_4(1)^{-1}. \\
  \end{array}
                                      $$
Due to Proposition $2.1$ we get:  $N\in \mathcal{M}_0(4\Lambda_2+\Lambda_4).$

\end{Example}

\begin{Theo} Let $\mathfrak{g}$ be of type $A$. Then, the map  
$$
\begin{array}{ccccccl}
&\kappa &   :   &\M&    \rightarrow  &\bigcup\limits_{\lambda\in P,s\in \mathbb{Z}} \mathcal{M}_s(\lambda)& \\
         &&          &M&       \mapsto    &(\Psi^{-1} \circ \Phi \circ \Psi)(M)&       
    \end{array}              
$$
defined as above is a morphism of crystals.
\end{Theo}

\begin{proof}
We have already seen that $\Psi$ is a crystal morphism therefore we limit the proof to $\Phi$ and get the claim by composition. Since $\Phi$ is successively defined it suffices to show that sending a reduced version matrix $M$ to $M^{(1)}$ preserves the crystal structure.   So we take such a matrix  $M=m_{ij}$ and its lower decomposition $M=M_1+M_2$ with $M_1=m_{ij}^{(1)}$ and $M_2=m_{ij}^{(2)}$. \\ 

By definition we obtain
$$(M^{(1)})_{ij}=m_{ij}^{(1)}+m_{i,j+1}^{(2)}.$$ 
Now we have to show the following for $i\in I$:

\begin{itemize}
\item[(i)]$wt(M)=wt(M^{(1)})$,
\item[(ii)]$\ph_i(M)=\ph_i(M^{(1)}),$
\item[(iii)]$\ep_i(M)=\ep_i(M^{(1)})$
\end{itemize}

and that computing $M^{(1)}$ interchanges with the Kashiwara operators, namely
\begin{center}
\begin{itemize}
\item[(iv)]$(\f_i(M))^{(1)}=\f_i(M^{(1)}),$
\item[(v)]$(\e_i(M))^{(1)}=\e_i(M^{(1)}).$
\end{itemize}
\end{center}
Since $\sum\limits_{j\in \Z} m_{ij}$ doesn`t change it is obvious that $wt$ is invariant under this construction. \\
We prove (ii) and (iv) simultaneously and (iii) and (v) follow in an analogous manner.   
For simplicity we set $M^{(1)}=:N=n_{ij}.$\\
We know that 
$$\ph_i(M)=\mbox{max} \{\sum\limits_{j \leq k}m_{ij} - \sum\limits_{j <k} m_{i+1,j} ;  k \in \Z\}.$$
Assume $\ph_i(M)\neq0$ and let $k\in \Z$ be minimal such that 
$$\ph_i(M)=\sum\limits_{j \leq k}m_{ij} - \sum\limits_{j <k} m_{i+1,j} .$$
Now assume $\ph_i(M)<\ph_i(N).$ \\
That is only possible if there exists $p\in\Z$ such that  $m_{i,p+1}^{(2)}>m_{i+1,p}^{(2)}$. \\
Otherwise $\sum\limits_{j \leq p}n_{ij} - \sum\limits_{j <p} n_{i+1,j}$ is equal to or  smaller than $\sum\limits_{j \leq p}m_{ij} - \sum\limits_{j <p} m_{i+1,j}.$\\ 
Due to lower decomposition this implies 
\begin{equation}
m_{i,p+1}^{(2)}\leq m_{i,p+1}-m_{i+1,p}^{(1)}
\end{equation}

and since $M=M_1+M_2,$
\begin{equation}
m_{i+1,p}^{(2)}= m_{i+1,p}-m_{i+1,p}^{(1)}.
\end{equation}
Now we compute 

$$
\begin{array}{ccl}
\sum\limits_{j \leq p}n_{ij} - \sum\limits_{j <p} n_{i+1,j} &=& \sum\limits_{j \leq p}m_{i,j}^{(1)}+ m_{i,j+1}^{(2)} - (\sum\limits_{j <p} m_{i+1,j}^{(1)} +m_{i+1,j+1}^{(2)}) \\
                                                                                              &=&\sum\limits_{j < p}m_{i,j}^{(1)}+ m_{i,j+1}^{(2)}+m_{i,p}^{(1)}+ m_{i,p+1}^{(2)} \\
                                                                                             && -
                                                                                             (\sum\limits_{j <p-1} m_{i+1,j}^{(1)} +m_{i+1,j+1}^{(2)}+ m_{i+1,p-1}^{(1)} +m_{i+1,p}^{(2)})\\
                                                                                              &\stackrel{(1)}{\leq}& 
                                                                                              \sum\limits_{j < p}m_{i,j}^{(1)}+ m_{i,j+1}^{(2)}+m_{i,p}^{(1)}+ m_{i,p+1}-m_{i+1,p}^{(1)} \\
                                                                                              &&- (\sum\limits_{j <p-1} m_{i+1,j}^{(1)} +m_{i+1,j+1}^{(2)}+ m_{i+1,p-1}^{(1)} +m_{i+1,p}^{(2)})\\
                                                                                              &\stackrel{(2)}{=}& \sum\limits_{j < p}m_{i,j}^{(1)}+ m_{i,j+1}^{(2)}+m_{i,p}^{(1)}+ m_{i,p+1}-m_{i+1,p}^{(1)} \\
                                                                                              &&- (\sum\limits_{j <p-1} m_{i+1,j}^{(1)} +m_{i+1,j+1}^{(2)} \\
                                                                                              &&+ m_{i+1,p-1}^{(1)} +m_{i+1,p}-m_{i+1,p}^{(1)})\\
                                                                                              &=& \sum\limits_{j \leq p+1}m_{ij} - \sum\limits_{j <p+1} m_{i+1,j}.

\end{array}
$$
But due to the choice of $k$ we have $$\ph_i(M)\geq\sum\limits_{j \leq p+1}m_{ij} - \sum\limits_{j <p+1} m_{i+1,j} $$
and hence
$$\ph_i(M)\geq\ph_i(N).$$
Now suppose $\ph_i(N)<\ph_i(M). $\\
It is obvious that $\sum\limits_{j \leq k}n_{ij} - \sum\limits_{j <k} n_{i+1,j}$ becomes smaller than $\sum\limits_{j \leq k}m_{ij} - \sum\limits_{j <k} m_{i+1,j}$ if and only if:
$$0\neq m_{i+1,k}^{(2)}>m_{i,k}^{(2)}.$$
Moreover the choice of $k$ again implies that
$$m_{i+1,k}^{(1)}\leq m_{i+1,k} <m_{i,k}.$$
That means we have the following situation in $M_1$: 
$$\sum\limits_{l\geq i+1}m_{l,k}^{(1)}= \sum\limits_{l > i+1}m_{l,k-1}^{(1)}$$
and therefore
$$m_{i,k}^{(2)}=m_{i,k}-m_{i+1,k-1}^{(1)}.$$
With this equation we compute:
$$
\begin{array}{ccl}
\sum\limits_{j \leq k-1}n_{ij} - \sum\limits_{j <k-1} n_{i+1,j} &=& \sum\limits_{j \leq k-1}m_{i,j}^{(1)}+ m_{i,j+1}^{(2)} - (\sum\limits_{j <k-1} m_{i+1,j}^{(1)} +m_{i+1,j+1}^{(2)}) \\
                                                                                              &=&\sum\limits_{j < k-2}m_{i,j}^{(1)}+ m_{i,j+1}^{(2)}+m_{i,k-1}^{(1)}+ m_{i,k}^{(2)} \\
                                                                                             && -
                                                                                             (\sum\limits_{j <k-2} m_{i+1,j}^{(1)} +m_{i+1,j+1}^{(2)}+  \\
                                                                                              &&m_{i+1,k-2}^{(1)} +m_{i+1,k-1}^{(2)})\\
                                                                                              &=& 
                                                                                              \sum\limits_{j < k-2}m_{i,j}^{(1)}+ m_{i,j+1}^{(2)}+m_{i,k-1}^{(1)}+ m_{i,k}-m_{i+1,k-1}^{(1)} \\
                                                                                              &&-(\sum\limits_{j <k-2} m_{i+1,j}^{(1)} +m_{i+1,j+1}^{(2)} + m_{i+1,k-2}^{(1)}\\
                                                                                              && +m_{i+1,k-1}-m_{i+1,k-1}^{(1)})\\
                                                                                             &=& \sum\limits_{j < k-2}m_{i,j}^{(1)}+ m_{i,j+1}^{(2)}+m_{i,k-1}^{(1)}+ m_{i,k} \\
                                                                                              &&-(\sum\limits_{j <k-2} m_{i+1,j}^{(1)} +m_{i+1,j+1}^{(2)} + m_{i+1,k-2}^{(1)} \\
                                                                                               && +m_{i+1,k-1})\\
                                                                                              &=& \sum\limits_{j \leq k}m_{ij} - \sum\limits_{j <k} m_{i+1,j}.

\end{array}
$$
Hence 
$$\ph_i(N)\geq\ph_i(M).$$
Furthermore these computations also show that we obtain $\ph_i(N)$ either by $$\sum\limits_{j \leq k}n_{ij} - \sum\limits_{j <k} n_{i+1,j}$$ or as in the last case by $$\sum\limits_{j \leq k-1}n_{ij} - \sum\limits_{j <k-1} n_{i+1,j} $$
and that $k$ (resp. $k-1$) is minimal with this property. \\

That means that we obtain $\f_i(N)$ by operating on $n_{i,k}$ or on $n_{i,k-1}.$\\
So assume we operate on $n_{i,k}$ and consider the lower decomposition of $\f_i(M)$. In this case we know that $m_{i,k}>m_{i+1,k-1}$ and the same in $N$. That  means $n_{i,k}>n_{i+1,k-1}.$ \\ 
Hence $$m_{i,k}^{(1)}>m_{i+1,k-1}^{(1)}.$$
 That implies $$\sum\limits_{j\geq i+2}m_{j,k-1}^{(1)}>\sum\limits_{j\geq i+1}m_{j,k}^{(1)}$$ 
and in particular $$m_{i+1,k}^{(2)} =0.$$
Therefore we get $m_{i+1,k}^{(1)}$ increased by one in the lower decomposition of $\f_i(M)$. \\
Moreover we know that 
$$m_{i+1,k}^{(1)}=m_{i+1,k} \geq m_{i,k+1}.$$ 
That means 
$$m_{i,k+1}^{(2)}=0$$ 
and increasing $m_{i+1,k}^{(1)}$ by one doesn`t change the decomposition of the $k+1$-st column.\\
Hence 
$$(\f_i(M))^{(1)}=\f_i(N).$$
Same arguments show that $m_{i+1,k}^{(2)}$ is  increased by one if we operate on $n_{i,k-1}$ and we also get   
$$(\f_i(M))^{(1)}=\f_i(N).$$
\qed
\end{proof}

For $M\in \M$ we set $B(M)$ to be the connected component of $M$ in $\M$. 
\begin{Cor} For $M \in \M$ and $\Phi(\Psi(M))=:n_{i,j} \in  \N$ we consider $s \in \mathbb{Z}$ maximal such that $n_{i,j}=0$ for all $j<s$ and $i\in \{1,\ldots,n+1\}$. Furthermore for $k=1,\ldots,n+1$ we define the values:
$$
a_k:=\sum\limits_{i\in \{1,\ldots,n+1\}} n_{i,k+s-1} - \sum\limits_{i\in \{1,\ldots,n+1\}} n_{i,k+s}\geq 0.
$$  
Then we have 
$$ \kappa(M) \in \mathcal{M}_s(\sum\limits_{k=1}^{n}a_k\Lambda_k)$$
and hence by restriction 

$$
\begin{array}{ccccccl}
&\kappa_{|B(M)} &   :   &B(M)&    \rightarrow  & \mathcal{M}_s(\sum\limits_{k=1}^{n}a_k\Lambda_k)& \\
    
     \end{array}              
$$
is a crystal isomorphism.

\end{Cor}

Now we give an application of the compression. In their framework about the correspondence between Young walls and Young tableaux, Kim and Shin [7] gave another realization of the crystal bases $B(\lambda)$ in the sense of reversed Young tableaux. Moreover Kang, Kim and Shin [2] constructed a crystal morphism between the monomials in $\mathcal{M}_1(\lambda)$ for dominant integral weights $\lambda$ and those reversed tableaux. By combining this with the crystal morphism $\kappa$ we can generalize their morphism to arbitrary monomials in $\M$.  \\
%
%
For a dominant integral weight $\lambda$ we define $S(\lambda)$ to be the set of all (reversed) semistandard tableaux of shape $\lambda$ with entries $1, \ldots , n+1$, which gives  a realization of  the crystal bases $B(\lambda) \ [7].$ \\ \\
Let $M\in \mathcal{M}_1(\lambda)$ be a monomial and $m_{ij}$ the associated reduced matrix in $\N$. We define the tableau $S(M)$ to be the semistandard reversed tableau with $m_{ij}$-many $i$ entries in $j$-th row. \\
Then we get 
\begin{Prop}$[2]$
The map
$$\begin{array}{ccccccl}
&\Omega&  :   &\mathcal{M}_1(\lambda)&  \rightarrow  &S(\lambda)& \\
&&                         &M&    \mapsto          &S(M)&
\end{array}
$$

is a crystal isomorphism.
\end{Prop}
It is obvious how to generalize this morphism to $\mathcal{M}_s(\lambda)$. Let $M\in \mathcal{M}_s(\lambda)$ be a monomial and $m_{ij}$ the associated reduced matrix in $\N$. In this case we define $S(M)$ to be the semistandard reversed tableaux with $m_{ij}$-many $i$ entries in the $j-s+1$-st row and get the morphism  
$$\begin{array}{ccccccl}
&\Omega&  :   &\bigcup\limits_{\lambda\in P,s\in \mathbb{Z}} \mathcal{M}_s(\lambda)&  \rightarrow  &\bigcup\limits_{\lambda \in P}S(\lambda)& \\
&&                         &M&    \mapsto          &S(M).&
\end{array}
$$

The combination of this morphism with the compression map $\kappa$ yields:
\begin{Cor}
The map
$$
\begin{array}{ccccccl}
&\Omega \circ \kappa  &  :   &\M&  \rightarrow  &\bigcup\limits_{\lambda \in P}S(\lambda)& \\
\end{array}
$$

is a crystal morphism.

\end{Cor}

\begin{Example} For $\mathfrak{g}$ of type $A_4$ we consider and the monomial $$M=Y_1(4)^{-1}Y_3(1) Y_1(3)^{-1}Y_4(1)^{-1}Y_2(0)^2 Y_3(2)^2.$$
We have already seen that $$\kappa (M)=N=X_1(1)^2X_1(3)X_2(0)^2X_2(1)X_3(0)^2X_3(1)X_5(0),$$ with the corresponding reduced matrix 
$$\begin{pmatrix}      
                                     0 & 2 & 0 &  1  \\ 
                                     2 & 2 & 1  & 0 \\       
                                     2 & 1 & 0  & 0 \\
                                     0 & 0 & 0  & 0  \\
                                     1 & 0 & 0  & 0  
       \end{pmatrix}.
$$
Therefore we assign the following semistandard reversed Young tableau
$$S(M)= \  {\young(::::1,::::2,11223,22335)} \ .$$

\end{Example}

\section{Insertion scheme for monomials in type $A$}

In this section, as another application of the compression given in section 2, we define a bumping rule for Nakajima monomials. That means we consider the crystal tensor product of two monomials $M_1$ and $M_2$ and search for a monomial $N \in \bigcup_{\lambda \in P,s \in \Z }\mathcal{M}_s(\lambda)$ such that the connected component of $M_1\otimes M_2$ is isomorphic to the connected component of $N$. Moreover we will see that this bumping is compatible with the reversed bumping for reversed tableaux (given in $[8]$).  \\
Before we define the monomial bumping we recall the tensor product rule for crystals $B_1$ and $B_2$: \\
The set $B_1 \otimes B_2 := \{b_1 \otimes b_2;b_1\in B_1$ and $b_2\in B_2\}$ becomes a crystal by setting

$$
\begin{array}{ccl}
 wt(b_1 \otimes b_2)          &=&    wt(b_1) + wt(b_2), \\ 
     \ep_i(b_1 \otimes b_2)       &=&    $max$\{\ep_i(b_1),\ep_i(b_2)+ \ll h_i,wt(b_1) \rr\}, \\
     \ph_i(b_1 \otimes b_2)       &=&    $max$\{\ph_i(b_1)+ \ll h_i,wt(b_2) \rr,\ph_i(b_2)\}, \\  \\
    
\f_i(b_1\otimes b_2) &=& \begin{cases}  \f_ib_1\otimes b_2 & \text{if  }  \ph_i(b_1)>\ep_i(b_2), \\
                                                             b_1\otimes \f_i  b_2 &  \text{if  }  \ph_i(b_1)\leq \ep_i(b_2) ,\\
                                  \end{cases}\\ \\

\e_i(b_1\otimes b_2) &=& \begin{cases}  b_1\otimes \e_ib_2 & \text{if  }  \ph_i(b_1)<\ep_i(b_2), \\
                                                             \e_ib_1\otimes   b_2 &  \text{if  }  \ph_i(b_1)\geq \ep_i(b_2). \\
                                  \end{cases}
\end{array}
$$

\vspace{2mm}
Let now $M_1$ and $M_2$ be reduced matrices of monomials in $\M$. In order to use the compression procedure we associate a matrix $M_1*M_2 \in  \mbox{Mat}_{n+1 \times \mathbb{Z}} (\Z_{\geq 0})$ to the tensor product $M_1\otimes M_2$ in the following way: 
$$
M_1\otimes M_2 \mapsto 
\begin{pmatrix}                          
& & & & 0 & & & & \\
& & & & 0 & & &  & \\
& &M_2 & & \vdots & &M_1 & & \\
& & & & 0 & & &  & \\
& & & & 0 & & &  &

\end{pmatrix}
=:M_1*M_2,
$$
where again $M_1$ and $M_2$ stand for their finite parts different from zero. \\
With the crystal structure on $\mbox{Mat}_{n+1 \times \mathbb{Z}} (\Z_{\geq 0})/ \sim$ we can show:

\begin{Prop} The map
$$
\begin{array}{ccccl} 
 &\mbox{Mat}_{n+1 \times \mathbb{Z}} (\Z_{\geq 0})/ \sim \otimes  \mbox{ Mat}_{n+1 \times \mathbb{Z}} (\Z_{\geq 0})/ \sim& \rightarrow &\mbox{Mat}_{n+1 \times \mathbb{Z}} (\Z_{\geq 0})/ \sim& \\
 
 &M_1\otimes M_2& \mapsto &M_1*M_2& 
 \end{array}
 $$ is a crystal morphism.
\end{Prop}
  
\begin{proof} Let $M_1$ and $M_2$ be matrices as above and after possible renumbering we set $M_1=(m_{i,j}^1)_{{i=1,\ldots,n+1}\atop {j=1,\ldots,l}}$ and $M_2=(m_{i,j}^2)_{{i=1,\ldots,n+1}\atop {j=1,\ldots,t}}.$ For simplicity we write \\ 
$M_1*M_2=:M=m_{i,j}$. \\
Since $M_1$ and $M_2$ are reduced we get by definition that $M_1*M_2$ is a reduced matrix without any (A1) application. \\
We have to show that $wt,\ph_i$ and $\ep_i$ are invariant under this map and that it commutes with $\f_i$ and $\e_i$. We observe directly from the definition that 
$$wt(M_1\otimes M_2)=wt(M_1*M_2).$$ 
For $i\in I$ we show that 
$$\ph_i(M_1\otimes M_2)=\ph_i(M_1*M_2).$$
In order to do this we distinguish the two cases $\ph_i(M_1)>\ep_i(M_2)$ and $\ph_i(M_1)\leq \ep_i(M).$
At first we assume $\ph_i(M_1)>\ep_i(M_2)$ and take $k$ minimal such that $$\ph_i(M_1)=\sum\limits_{j\leq k} m_{i,j}^1-m_{i+1,j-1}^1.$$
This implies 
$$
\begin{array}{ccccl}
&\ph_i(M_1*M_2)&                                                                &=& \sum\limits_{j\leq k+t+1} m_{i,j}-m_{i+1,j-1}. \\ 
                          &&                                                                  &=& \sum\limits_{j\leq k} m_{i,j}^1-m_{i+1,j-1}^1 \\
                            &&                                                                     &&+ \sum\limits_{j} m_{i,j}^2 -\sum\limits_{j} m_{i+1,j}^2 \\                                                                                                     &\mbox{since } \ll h_i,\Lambda_j \rr=\delta_{i,j}, &                   &=& \ph_i(M_1) +\sum\limits_{j} m_{i,j}^2 \ll h_i, \Lambda_i-\Lambda_{i-1} \rr \\
                           &&                                                                  &&    + \sum\limits_{j} m_{i+1,j}^2 \ll h_i, \Lambda_{i+1}-\Lambda_{i} \rr  \\
                           &&                                                                  &=& \ph_i(M_1) + \sum\limits_{j} m_{i,j}^2 \ll h_i, \beta_i \rr \\
                            &&                                                                        &&+ \sum\limits_{j} m_{i+1,j}^2 \ll h_i, \beta_{i+1} \rr  \\
                           &&                                                                  &=&  \ph_i(M_1)  \\
                           &&                                                                      && +\ll h_i, \sum\limits_{j} m_{i,j}^2 \ \beta_i + \sum\limits_{j} m_{i+1,j}^2  \ \beta_{i+1}\rr   \\ 
& \ll h_i,\beta_j \rr=0 \mbox{ for } j\neq i,i+1,&      &=& \ph_i(M_1) + \ll h_i, \sum\limits_{i} (\sum\limits_{j} m_{i,j}^2)\beta_i \rr \\
                      &&                                                                        &=& \ph_i(M_1) + \ll h_i,wt(M_2) \rr \\
                   &\mbox{ since }  \ph_i(M_1)>\ep_i(M_2),&              &=& \ph_i (M_1\otimes M_2).
\end{array}
$$
Moreover these computations also show: If $\ph_i(M_1)>\ep_i(M_2),$ then $k$ is minimal with
$$\ph_i(M_1*M_2)=\ph_i(M)=\sum\limits_{j\leq k} m_{i,j}-\sum\limits_{j<k}m_{i+1,j}$$

and hence 
$$\f_i(M_1*M_2)=\f_i M_1*M_2.$$ 
Due to the tensor product rule we also observe 
$$\f_i(M_1 \otimes M_2)=\f_iM_1 \otimes M_2,$$
which implies that the map interchanges with $\f_i$ in this case. \\

Let now $\ph_i(M_1)\leq\ep_i(M_2)$. This yields directly $$\ph_i(M_1*M_2)=\ph_i(M_2)=\ph_i(M_1 \otimes M_2)$$ 
and again $\f_i(M_1*M_2)=M_1*\f_i M_2.$  \\
Same arguments hold for $\ep_i$ and $\e_i$.

\qed 
\end{proof}

With this interpretation of the crystal tensor product of monomials we are able to give the definition of bumping for Nakajima monomials. Let $M_1$ and $M_2$ be monomials in $\M$ then we define $M_1 \rightarrow M_2$ as the result of the following compositions of crystal morphisms:
$$
\begin{array}{ccccccccl}
&\M \times \M& \rightarrow     &\mbox{Mat}_{n+1 \times \mathbb{Z}} (\Z_{\geq 0})/ \sim&  \stackrel{\Phi}{\rightarrow}      &\N/ \sim&         \rightarrow    &\bigcup\limits_{\lambda\in P,s\in\mathbb{Z}}\mathcal{M}_s(\lambda)&                                             \\
 &(M_1,M_2)&   \mapsto        &\Psi(M_1)*\Psi(M_2)&                                                          \stackrel{\Phi}{\mapsto}           &N& \mapsto   &\Psi^{-1}(N).&
\end{array}
$$
In other words we set
$$M_1 \rightarrow M_2:=\Psi^{-1}(\Phi(\Psi(M_1)*\Psi(M_2))).$$
With Proposition 3.1 and Theorem 2.1 we observe

\begin{Theo}Let $\mathfrak{g}$ be of type $A$. Then, the map  
$$
\begin{array}{ccccccl}
     &\M \otimes \M&    \rightarrow  &\bigcup\limits_{\lambda\in P,s\in \mathbb{Z}} \mathcal{M}_s(\lambda)& \\
       &M_1 \otimes M_2&       \mapsto    &M_1\rightarrow M_2&       
    \end{array}              
$$
defined as above is a morphism of crystals.

\end{Theo}
Now we notice that the monomial bumping coincides with the tableaux bumping defined in $[8]$. More precisely, if we take $M_1,M_2 \in \M$ we have two possibilities to associate a reversed tableaux to their tensor product $M_1\otimes M_2$. \\
The first one is to take the monomial bumping $M_1\rightarrow M_2$ and to consider the tableaux $S(M_1\rightarrow M_2)$ in the sense of Corollary 2.2. On the other hand we compute $S(M_1)$ and $S(M_2)$ and apply the reversed bumping rule given by Kim and Shin [7], [8], namely $S(M_1) \rightarrow S(M_2)$.  Corollary 2.2  and Theorem 3.1 imply
$$  S(M_1\rightarrow M_2)= S(M_1) \rightarrow S(M_2). $$

\section{Compression of Nakajima monomials in type $C$}
In this section we will define the compression of Nakajima monomials for a Lie algebra $\frak{g}$ of type $C_n$. We briefly recall the basic setting of $\frak{g}$. Let $P$ be the weight lattice of $\frak{g}$ and $\beta_1,\ldots,\beta_n$ the orthogonal basis of $P$. Let further $I=\{1,\ldots,n\}$ be the index set for the simple roots given by $\alpha_i=\beta_i-\beta_{i+1}$ for $i=1,\ldots,n-1$ and $\alpha_n=2\beta_n$. Moreover we get the fundamental weights by $\Lambda_i=\beta_1+\ldots+\beta_i$ and therefore $\beta_i=\Lambda_i-\Lambda_{i-1}$. Then we compute for all $i\neq n$ and $j\in I$:
$$ \left\langle h_j, \alpha_i\right\rangle=
\begin{cases}
  2  &\text{if $i=j$,} \\
  -1 &\text{if $j=i-1$ or $j=i+1$,} \\
  0  &\text{else}
\end{cases}
$$
and
$$\left\langle h_j, \alpha_n\right\rangle=
\begin{cases}
  2  &\text{if $j=n$,} \\
  -2 &\text{if $j=n-1$,} \\
  0  &\text{else.}
\end{cases}
$$
As in the $A_n$-case we set 
$$c_{ij}= \begin{cases}     0    &\text{if  $i>j$,} \\      
                                       1    &\text{else.}  
              \end{cases}                         
$$               
With this notation we obtain for $j\in\Z$
$$ 
A_i(j)= \begin{cases}   Y_i(j)Y_{i+1}(j)^{-1}Y_i(j+1)Y_{i-1}(j+1)^{-1}      &\text{if  $i \neq n$,} \\      
                                       Y_n(j)Y_n(j+1)Y_{n-1}(j+1)^{-2}           &              \text{if  $i = n$.} 
              \end{cases}                         
$$               
Let $\B=\{1,\ldots,n,\bar{1},\ldots\bar{n}\}$ then we define a total order on $\B$ by 
$$1\prec2\prec\ldots\prec n\prec \overline{n} \prec \ldots \prec \overline{2} \prec \overline{1}.$$
For $i\in I$ and $j\in \Z$ we consider the variables defined in $[3]$:
$$
\begin{array}{ccl}
   X_i(j)&:=&Y_{i-1}(j+1)^{-1}Y_i(j), \\
   X_{\overline{i}}(j)&:=& Y_{i-1}(j+n-i+1)Y_i(j+n-i+1)^{-1}.
\end{array}
$$   
With these variables we have for $i\neq n:$
$$
\begin{array}{ccl}
   A_i(j)  &=&  X_{i+1}(j)^{-1}X_i(j), \\
   A_i(j)  &=&  X_{\overline{i+1}}(j-n+i)X_{\overline{i}}(j-n+i)^{-1} \\
\end{array}
$$
and  
$$
A_n(j)    =   X_n(j)X_{\overline{n}}(j)^{-1}.
$$
Furthermore it is easy to see that for $i=1,\ldots,n$ and some $p-q=n-i$ the following equation holds
$$ X_i(p)X_{\overline{i}}(q)=X_{i+1}(p)X_{\overline{i+1}}(q).$$
This equation will be important later when we define the equivalence relation on matrices because it involves more options to write an arbitrary monomial as a product of  $X_i(j)$`s and $X_{\overline{i}}(j)$`s.  \\
As in section 2 we recall the characterization of $M_1(\lambda)$ with a dominant integral weight $\lambda$ for Lie algebras of type $C_n$ given in $[3]$:

\begin{Prop}$[3]$
Let $\lambda=a_1\Lambda_1+ \ldots +a_n\Lambda_n$. Then the connected com-ponent $M_1(\lambda)$ containing the maximal vector
$$M_1=Y_1(1)^{a_1}\cdots Y_n(1)^{a_n}$$
 is characterized as the set of monomials 
 $$M=X_{t_{1,1}}(1)\cdots X_{t_{1,k_1}}(1)\cdots X_{t_{n,1}}(n)\cdots X_{t_{n,k_n}}(n)$$
 satisfying the following conditions:
 \begin{itemize}
 \item[(i)] $k_j=a_j+\ldots+a_n$ for all $j=1,\ldots,n,$
 \item[(ii)]$t_{j,1} \succeq t_{j,2} \succeq \ldots \succeq t_{j,k_j}$ for all $j=1,\ldots,n,$ 
 \item[(iii)] for each $j=2,\ldots,n$ and $l=1,\ldots,k_j, \  t_{j-1,l}\succ t_{j,l}.$
 
 \end{itemize}
 \end{Prop}
 For $s\in \mathbb{Z}$ we also consider the \emph{shifted} highest weight monomials of weight $\lambda$ $$M_s=Y_1(s)^{a_1}Y_2(s)^{a_2} \ldots Y_n(s)^{a_n}.$$ 
As an immediate consequence of Proposition 4.1 we obtain their connected component $\mathcal{M}_s(\lambda)$ by the set of monomials of the form 
$$M=X_{t_{1,1}}(s)\cdots X_{t_{1,k_1}}(s)\cdots X_{t_{n,1}}(s+n-1)\cdots X_{t_{n,k_n}}(s+n-1)$$
 
satisfying condition $(i),(ii)$ and $(iii).$ \\

We will see later on that these conditions translate into the notation of matrices exactly the same way as in the $A_n$-case.
In order to use similar constructions as in section 2 we show that $\M$ is generated by the elements $X_i(j)$ and $X_{\overline{i}}(j)$. That means we define $\mathbf{M}$ to be the free abelian monoid  generated by the set \{$X_i(j),X_{\overline{i}}(j),i=1,\ldots,n,j\in \Z \}.$ We further define an ideal $ \mathbf{J} \subset \mathbf{M}$ by 
$$
\begin{array}{ccl}
\mathbf{J} &:=& \langle \prod\limits_{k=1}^{i} X_k(j+i-k) X_{\overline{k}}(j-n+k-1), \\
                  && \prod\limits_{k=1}^{i} X_{\overline{k}}(j-i+k) X_k(j+n+1-k),i =1,\ldots,n, j\in \Z \rangle_{\mathbf{M}}.
\end{array}
$$
Hence the quotient $\mathbf{M}/\mathbf{J}$ becomes a group with
$$
X_i(j)^{-1}=\prod\limits_{k=1}^{i-1} X_k(j+i-k) \prod\limits_{k=1}^{i} X_{\overline{k}}(j-n+k-1)
$$
and 
$$
X_{\overline{i}}(j)^{-1}=\prod\limits_{k=1}^{i-1} X_{\overline{k}}(j-i+k)\prod\limits_{k=1}^{i} X_k(j+n+1-k).
$$
This gives rise to an analog of Proposition 2.2.

\begin{Prop} We have 
$$
\M \cong \mathbf{M}/\mathbf{J}.$$
\end{Prop}

\begin{proof} We consider the map that identifies $X_i(j)$ with $Y_i(j)Y_{i-1}(j+1)^{-1}$ and $X_{\overline{i}}(j)$ with $Y_{i-1}(j+n-i+1)Y_i(j+n-i+1)^{-1}$.
Let $M=\prod\limits_{i\in I,j \in \Z} Y_i(j)^{y_i(j)}$ be a monomial in $\M$. In order to show surjectivity we consider again each $Y_i(j)^{y_i(j)}$ separately and distinguish two cases:
 
\underline{1. case}: $y_i(j)>0$, then we write 

$$
\begin{array}{ccl}
Y_i(j)&=&\prod\limits_{k=0}^{i-1} Y_k(j+i-k)^{-1} \prod\limits_{k=1}^{i}Y_k(j+i-k)\\ &=& \prod\limits_{k=1}^{i} X_k(j+i-k).
\end{array}
$$

Therefore we get

$$Y_i(j)^{y_i(j)}= \prod\limits_{k=1}^{i} X_k(j+i-k)^{y_i(j)}.$$

\underline{2. case}: $y_i(j)<0$, then we set

$$
Y_i(j)^{-1} = \prod\limits_{k=1}^{i} X_{\overline{k}}(j-n+k-1)         
$$

and hence 
$$
Y_i(j)^{y_i(j)}= \prod\limits_{k=1}^{i} X_{\overline{k}}(j-n+k-1)^{-y_i(j)}. 
$$
With these equations we compute 
$$
X_i(j)^{-1}=Y_i(j)^{-1}Y_{i-1}(j+1)= \prod\limits_{k=1}^{i-1} X_k(j+i-k) \prod\limits_{k=1}^{i} X_{\overline{k}}(j-n+k-1)
$$
and
$$
\begin{array}{ccl}
X_{\overline{i}}(j)^{-1}&=&Y_{i-1}(j+n-i+1)^{-1}Y_i(j+n-i+1)\\
                                 &=& \prod\limits_{k=1}^{i-1} X_{\overline{k}}(j-i+k)\prod\limits_{k=1}^{i} X_k(j+n+1-k),
\end{array}
$$
which implies injectivity. 

\qed
\end{proof}
\begin{Remark}
The equation 
$$ X_i(p)X_{\overline{i}}(q)=X_{i+1}(p)X_{\overline{i+1}}(q)$$
for $p-q=n-i$, also holds in $\mathbf{M}/\mathbf{J}$.
\end{Remark}


Due to Proposition 4.2 we can write an arbitrary $M\in \M$ as a product of $X_i(j)$`s and $X_{\overline{i}}(j)$`s. More precisely, there exist $m_{ij} \in \Z_{\geq 0}$ such that 
 $$M=\prod\limits_{i\in \B,j \in \Z} X_i(j)^{m_{ij}}.$$ 

In other words, we can associate a matrix $m_{i,j}$ to each $M\in \M$, where $i\in \{1,\ldots,n,\bar{1},\ldots\bar{n}\}$ and $j\in \Z.$ That means we obtain a matrix in  $\mbox{Mat}_{2n \times \Z} (\Z_{\geq 0}).$ As in the $A_n$-case these are matrices with just finitely many non zero columns and we number the rows by $1,\ldots,n,\bar{n},\ldots\bar{1}$ instead of $1,\ldots,2n$. \\
In order to get a bijection between the monomials and those matrices we need to fix the matrix notation. Consider  $M\in \mbox{Mat}_{2n \times \Z} (\Z_{\geq 0})$ with $M=m_{i,j}$. Then the definition of $\mathbf{J}$ and Remark 4.1 allow us to apply the following rules without changing the underlying monomial:
\begin{itemize}

 \item[\textbf{(C1)}] For a pair $p,q$ with $p-q=n-b$ with
 $$ m_{b,p}\neq 0 \mbox{  and  }  m_{\overline{b},q}\neq 0$$
 we decrease $ m_{b,p}$ and $m_{\overline{b},q}$ by min$\{m_{b,p}, m_{\overline{b},q}\}$
 and increase $ m_{b+1,p}$ and  $m_{\overline{b+1},q}$ by min$\{m_{b,p}, m_{\overline{b},q}\}.$

\item[\textbf{(C2)}] 
For a pair $p,q$ with $p-q=n-b+1$ with
 $$ m_{b,p}\neq 0 \mbox{  and  }  m_{\overline{b},q}\neq 0$$
 we decrease $ m_{b,p}$ and $m_{\overline{b},q}$ by min$\{m_{b,p}, m_{\overline{b},q}\}$
 and increase $ m_{b-1,p}$ and  $m_{\overline{b-1},q}$ by min$\{m_{b,p}, m_{\overline{b},q}\}.$
\end{itemize}
\vspace{3mm}
Moreover we have an analog of rule (A1):
\vspace{3mm}
\begin{itemize}

 \item[\textbf{(C3)}] 
For every $i\in \{ 1, \ldots,n \}$ and $k\in \Z$ with 
\begin{center}
 $m_{i-s,k+s} \neq 0 $ for all $s=0,1,\ldots, i-1$
 \end{center} 
 and 
\begin{center} 
 $m_{\overline{i-s},k-n+i-1-s} \neq 0 $ for all $s=0,1,\ldots, i-1$
\end{center}
we decrease all these entries by 
 $$
 \mbox{min}\{m_{i-s,k+s},m_{\overline{i-s},k-n+i-1-s}; 0,1,\ldots, i-1 \}.
 $$
 We call such a collection a generalized diagonal at $m_{i,k}$ and this procedure a cancellation at $m_{i,k}$. 
 
 \end{itemize}
 \vspace{3mm}
We can also insert some generalized diagonals to get longer ones:
\vspace{3mm}
\begin{itemize}

 \item[\textbf{(C4)}] 
For every $i\in \{ 1, \ldots,n \}$ and $k\in \Z$ with 
\begin{center}

$m_{i,k} \neq 0$ 
 \end{center} 
 and 
\begin{center} 
 $m_{\overline{i-s},k-n+i-1-s} \neq 0 $ for all $s=0,1,\ldots, i-1$
\end{center}
we increase the entries $m_{i-s,k+s}$ for all $s=1,\ldots, i-1$ by
$$
\mbox{min}\{m_{i,k},m_{\overline{i-s},k-n+i-1-s}; 0,1,\ldots, i-1 \}
 $$
 and apply (C3) to get a longer cancellation at $m_{i,k}$.    
 
\vspace{3mm}
For every $i\in \{ 1, \ldots,n \}$ and $k\in \Z$ with
$$m_{\overline{i},k-n+i-1} \neq 0$$
and 
\begin{center} 
 $m_{i-s,k+s} \neq 0$ for all $s=0,1,\ldots, i-1$
\end{center}
we increase the entries $m_{\overline{i-s},k-n+i-1-s}$ for all $s=1,\ldots, i-1$ by
$$
\mbox{min}\{m_{i-s,k+s},m_{\overline{i},k-n+i-1}; 0,1,\ldots, i-1 \}
 $$
 and apply (C3) to get a longer cancellation at $m_{i,k}$.   

\end{itemize}

\vspace{3mm}
We use the rules (C1)-(C4) to get  \emph{reduced versions} of the matrices associated to monomials in $\M$. We explain what we mean by reduced in this case:\\
Let  $m_{i,j}\in \mbox{Mat}_{2n \times \Z} (\Z_{\geq 0})$ be a matrix corresponding to a monomial $M\in \M$:
$$M=\prod\limits_{i\in \B,j \in \Z} X_i(j)^{m_{ij}}.$$
We search for $[m_{i,j}]$ such that 
\begin{itemize}
\item[(i)]
$\sum\limits_{i\in \B,j \in \Z}[m_{ij}]=\mbox{min} \{\sum\limits_{i\in \B,j \in \Z}n_{ij} \ ; \ \prod\limits_{i\in \B,j \in \Z} X_i(j)^{m_{ij}}=\prod\limits_{i\in \B,j \in \Z} X_i(j)^{n_{ij}}\}, $
\item[(ii)] there are no pairs $p,q$ with $p-q=n-b+1$ such that
 $$ [m_{b,p}]\neq 0 \mbox{  and  }  [m_{\overline{b},q}]\neq 0.$$
 \end{itemize}   
We use the rules (C1)-(C4) stepwise to obtain $[m_{ij}]$ from $m_{ij}$ as follows: 

\vspace{3mm}

Let $j$ be minimal such that $m_{i,k}=0$ for all $i \in \B$ and $k>j$. Then we start at $m_{n,j}$ and apply (C1)-(C2) to all the other entries if this yields an application of (C3)-(C4) to $m_{n,j}$. This means we try to get some cancellation at this entry. \\
After that we go left to the next entry in this row and do the same. \\
Once we have done this with the whole row we go to the upper one and apply the same procedure until we reach $m_{1,1}$.  \\
At the end we apply (C2) to guarantee the desired condition (ii).

With this notation we define an equivalence relation on $\mbox{Mat}_{2n \times \Z} (\Z_{\geq 0})$:
$$m_{i,j} \sim n_{i,j} \mbox { iff }  [m_{i,j}] = [n_{i,j}].$$

We consider the quotient  
\begin{center}
$\mbox{\mbox{Mat}}_{2n \times \mathbb{Z}} (\Z_{\geq 0})/ \sim$ .
\end{center}
It is obvious that two matrices which lie in the same equivalence class correspond to the same monomial and vice versa. Hence by sending a monomial onto its reduced matrix we get


\begin{Prop}
There exists a bijection between $\M$ and $\text{Mat}_{2n \times \mathbb{Z}} (\Z_{\geq 0})/ \sim$.   
\end{Prop}

In order to get a morphism of crystals we endow $\mbox{Mat}_{2n \times \mathbb{Z}} (\Z_{\geq 0})/ \sim$ with a crystal structure by defining it on the reduced representatives.  \\ \\
Let $M=m_{i,j} \in \mbox{Mat}_{2n \times \mathbb{Z}} (\Z_{\geq 0})$ be a reduced matrix. Then we set
$$
wt(M) = \sum\limits_{i=1}^{n}\big(\sum\limits_{j\in\Z}(m_{ij}-m_{\overline{i},j})\big)\beta_i. \\ 
$$
For $i\neq n$ we put
$$
\begin{array}{ccl}
\ph_i(M) &=&  $ max $ \{\sum\limits_{j \leq k}m_{ij}+m_{\overline{i+1},j-n+i} - \sum\limits_{j <k} m_{i+1,j}+m_{\overline{i},j-n+i}$ ; $ k \in \Z\}, \\      
 \ep_i(M) &=& -$ min $\{ \sum\limits_{j > k}m_{ij}+m_{\overline{i+1},j-n+i}- \sum\limits_{j \geq k} m_{i+1,j}+m_{\overline{i},j-n+i} $ ; $ k \in \Z\} \\
\end{array}
$$
and 
$$
\begin{array}{ccl}
\ph_n(M) &=&  $ max $ \{\sum\limits_{j \leq k}m_{n,j}- \sum\limits_{j <k}m_{\overline{n},j}$ ; $ k \in \Z\}, \\      
 \ep_n(M) &=& -$ min $\{ \sum\limits_{j > k}m_{n,j}- \sum\limits_{j \geq k} m_{\overline{n},j} $ ; $ k \in \Z\}. \\
\end{array}
$$
If $\ph_i(M)=0$ we set $\f_i(M)=0$ for all $i\in I$. \\
Let now $\ph_i(M)\neq0$ then we define the Kashiwara operator $\f_i$ for $i\neq n$: \\
Let $k$ be minimal such that 
$$
\ph_i(M) = \sum\limits_{j \leq k}m_{ij}+m_{\overline{i+1},j-n+i} - \sum\limits_{j <k} m_{i+1,j}+m_{\overline{i},j-n+i}.
$$
Then we distinguish the following two cases: \\ \\
\underline{1. case}:
$m_{\overline{i+1},k-n+i}=0.$\\
Then we set $\f_i(M)$ as the matrix we get from $M$ by increasing (resp. decreasing) $m_{i+1,k}$ (resp. $m_{i,k}$) by one. 
Formally spoken we obtain $\f_i(M)=\hat{m}_{s,j} $  from $M=m_{s,j}$ by 
\begin{center}
            $ \hat{m}_{s,j}  = \begin{cases} 
                                                   m_{s,j}            & \text{if  $(s,j) \notin \{(i,k),(i+1,k)\}$} , \\
                                               m_{i,k}-1              & \text{if  $(s,j)=(i,k)$,} \\
                                               m_{i+1,k}+1           & \text{if $(s,j)=(i+1,k)$.} \\
                                        \end{cases}    
            $
\end{center}                                            

\underline{2. case}: $m_{\overline{i+1},k-n+i}\neq0.$ \\
Then we define $\f_i(M)$ as the matrix we get from $M$ by increasing (resp. decreasing) $m_{\overline{i},k-n+i}$ (resp. $m_{\overline{i+1},k-n+i}$) by one. That means we obtain $\f_i(M)=\hat{m}_{s,j} $  from $M=m_{s,j}$ by 
\begin{center}
            $ \hat{m}_{s,j}  = \begin{cases} 
                                                   m_{s,j}            & \text{if  $(s,j) \notin \{(\overline{i},k-n+i),(\overline{i+1},k-n+i)\}$} , \\
                                               m_{\overline{i+1},k-n+i}-1              & \text{if  $(s,j)=(\overline{i+1},k-n+i)$,} \\
                                               m_{\overline{i},k-n+i}+1           & \text{if $(s,j)=(\overline{i},k-n+i)$.} \\
                                            \end{cases}    
            $
\end{center} 
Now we give the definition of $\f_n(M)$ for $\ph_n(M) \neq 0$: \\
Let $k$ be minimal such that
$$\ph_n(M) = \sum\limits_{j \leq k}m_{n,j}- \sum\limits_{j <k}m_{\overline{n},j}.$$
Then we set $\f_n(M)$ as the matrix we get from $M$ by increasing (resp. decreasing) $m_{\overline{n},k}$ (resp. $m_{n,k}$) by one.   
More precisely we obtain $\f_n(M)=\hat{m}_{s,j} $  from $M=m_{s,j}$ by \\
\begin{center}
            $ \hat{m}_{s,j}  = \begin{cases} 
                                                   m_{s,j}            & \text{if  $(s,j) \notin \{(n,k),(\overline{n},k)\}$} , \\
                                               m_{n,k}-1              & \text{if  $(s,j)=(n,k)$,} \\
                                               m_{\overline{n},k}+1           & \text{if $(s,j)=(\overline{n},k)$.} \\
                                            \end{cases}    
            $
\end{center}  
If $\ep_i(M)=0$ we set $\e_i(M)=0$.  \\    
For $\ep_i(M) \neq 0$ let $p$ be maximal such that 
$$\ep_i(M) =- (\sum\limits_{j > p}m_{ij}+m_{\overline{i+1},j-n+i}- \sum\limits_{j \geq p} m_{i+1,j}+m_{\overline{i},j-n+i}). $$
Then we distinguish the following two cases to define $\e_i(M)$ for $i \neq n$: \\ \\
\underline{1. case}: $m_{\overline{i+1},p-n+i} \neq 0.$\\
Then we set $\e_i(M)$ as the matrix we get from $M$ by increasing (resp. decreasing) $m_{i,p}$ (resp. $m_{i+1,p}$) by one. Formally spoken we observe $\e_i(M)=\hat{m}_{s,j} $  from $M=m_{s,j}$ by \\
\begin{center}
            $ \hat{m}_{s,j}  = \begin{cases} 
                                                   m_{s,j}            & \text{if  $(s,j) \notin \{(i,p),(i+1,p)\}$} , \\
                                               m_{i,p}+1              & \text{if  $(s,j)=(i,p)$,} \\
                                               m_{i+1,p}-1           & \text{if $(s,j)=(i+1,p)$.} \\
                                            \end{cases}    
            $
\end{center}                                            
\underline{2. case}: $m_{\overline{i+1},p-n+i} = 0.$\\
Then we define $\e_i(M)$ as the matrix we get from $M$ by increasing (resp. decreasing) $m_{\overline{i+1},p-n+i}$ (resp. $m_{\overline{i},p-n+i}$) by one. That means  we obtain $\e_i(M)=\hat{m}_{s,j} $  from $M=m_{s,j}$ by \\
\begin{center}
            $ \hat{m}_{s,j}  = \begin{cases} 
                                                   m_{s,j}            & \text{if  $(s,j) \notin \{(\overline{i},p-n+i),(\overline{i+1},p-n+i)\}$} , \\
                                               m_{\overline{i+1},p-n+i}+1              & \text{if  $(s,j)=(\overline{i+1},p-n+i)$,} \\
                                               m_{\overline{i},p-n+i}-1           & \text{if $(s,j)=(\overline{i},p-n+i)$.} \\
                                            \end{cases}    
            $
\end{center}                                            
Let $p$ be maximal such that
$$\ep_n(M) = -(\sum\limits_{j > p}m_{n,j}- \sum\limits_{j \geq p} m_{\overline{n},j}). $$                                      
Then we set $\e_n(M)$ as the matrix we get from $M$ by decreasing (resp. increasing) $m_{\overline{n},p}$ (resp. $m_{n,p}$) by one.  
Formally spoken we obtain $\e_n(M)=\hat{m}_{s,j} $  from $M=m_{s,j}$ by \\
\begin{center}
            $ \hat{m}_{s,j}  = \begin{cases} 
                                                   m_{s,j}            & \text{if  $(s,j) \notin \{(n,p),(\overline{n},p)\}$} , \\
                                               m_{n,p}+1              & \text{if  $(s,j)=(n,p)$,} \\
                                               m_{\overline{n},p}-1           & \text{if $(s,j)=(\overline{n},p)$.} \\
                                            \end{cases}    
            $
\end{center}       
Easy computations show that $\mbox{Mat}_{2n \times \mathbb{Z}} (\Z_{\geq 0})/ \sim$ along with the maps $wt,\ph_i,\ep_i,\f_i$ and $\e_i$ becomes a semi-normal crystal.
As in section 2 we prove that this crystal structure coincides with the structure on $\M$ under the above bijection.
\begin{Prop} The bijection
\begin{center}
$\begin{array}{ccccl}
 \Psi:  &\M&    &\rightarrow&  \mbox{Mat}_{2n \times \mathbb{Z}} (\Z_{\geq 0})/ \sim   \\
&M=\prod\limits_{i\in \B,j \in \Z} X_i(j)^{m_{ij}}&   &\mapsto&    [m_{ij}]
\end{array}$
\end{center}
is a crystal isomorphism. 
\end{Prop}

\begin{proof} It is easy to verify that $wt,\ep_i$ and $\ph_i$ are invariant under $\Psi$ especially because they are invariant under the application of (C1)-(C4). It remains to show that $\Psi $ commutes with the crystal operators $\f_i$ and $\e_i$. Let $M \in \M$ be a monomial with associated matrix $[m_{l,j}]$. Due to the crystal structure defined on the matrices it follows almost directly that for all $i\in I$
$$\f_i(M)=\prod\limits_{l\in \B,j \in \Z} X_l(j)^{\f_i([m_{lj}])}$$
and the same for $\e_i$.  \\
Therefore it suffices to verify that 
$$ \f_i([m_{lj}])=[\f_i([m_{lj}])]$$ 
and
$$ \e_i([m_{lj}])=[\e_i([m_{lj}])].$$ 
Since this can be proved analogously we just give the proof for $\f_i$.\\
For simplicity we denote $[m_{l,j}]$ by $m_{l,j}$ and for $i\neq n$ let $k$ be minimal such that 
$$
\ph_i(M) = \sum\limits_{j \leq k}m_{ij}+m_{\overline{i+1},j-n+i} - \sum\limits_{j <k} m_{i+1,j}+m_{\overline{i},j-n+i}.
$$ 
First we look at the case $m_{\overline{i+1},k-n+i}=0.$ That means we get $\f_i(m_{l,j})$ by increasing $m_{i+1,k}$ and decreasing $m_{i,k}$ each by one. Let us assume that the  increase of  $m_{i+1,k}$ induces a longer or new cancellation. But this yields 
$$m_{i+1,k}=0 \mbox{ and }  m_{\overline{i},k-n+i}=0$$ 

and together with the choice of $k$ this implies
$$m_{i,k+1}=0 \mbox{ and } m_{\overline{i+1},k-n+i+1}=0.$$
Since $m_{i,k+1}=0$  and $m_{\overline{i},k-n+i}=0$, the increase doesn`t provide any generalized diagonal at $m_{i+s,k+1-s}$ for positive $s$, without insertion. But those insertions would have been done before we operate with $\f_i$ because $m_{l,j}$ is in reduced form. \\
It is still possible that we get a new generalized diagonal by (C4) at $m_{i+1,k}$  itself. But we get no diagonal $m_{i+1-s,k+s}\neq 0$ for $s=0,\ldots, i$ since $m_{i,k+1}=0$. Furthermore we can`t increase this entry by applying (C1) because otherwise we could have applied (C4) at $m_{i,k+1}$ before. \\    
The second possibility to apply (C4) at $m_{i+1,k}$ needs $m_{\overline{i+1-s},k-n+i-s}\neq 0$ for all $s=0,\ldots, i$. In particular we get $m_{\overline{i-s},k-n-1+i-s}\neq 0$ for all $s=0,\ldots, i-1$ and $m_{i,k}\neq0$ which implies an application of (C4) at $m_{i,k}$ before operating. This is again a contradiction to the fact that $m_{i,j}$ is reduced. \\
Moreover since $m_{\overline{i+1},k-n+i}=0$ we can not apply (C1) to the increased $m_{i+1,k}$.\\
Overall we have seen that operating with $\f_i$ preserves the reduced version in this case. \\ 
Similar arguments hold for the case $m_{\overline{i+1},k-n+i} \neq 0$ .  

\qed
\end{proof}

Now we translate the characterizing conditions of the monomials that give a realization of the crystal bases, given in Proposition 4.1, into the language of matrices. We will recognize that those are the same conditions as in section 2. From this observation one can deduce  that similar constructions yield our desired morphism. 

\begin{Def}
Define $\N\subset \mbox{Mat}_{2n \times \mathbb{Z}} (\Z_{\geq 0}) $ as the set of matrices whose reduced versions have only zero-entries out of a $ 2n\times n$-submatrix \\ $M=(m_{ij})_{{i=1,\ldots,n,\overline{n},\ldots,\overline{1}}\atop{j=0,\ldots,n-1}}$ satisfying the following properties:
\begin{itemize}
\item[(i)]  $m_{ij} \in \Z_{ \geq 0}$ for $i=1,\ldots,n,\overline{n},\ldots,\overline{1}$ and $j=0,\ldots,n-1$,
\item[(ii)] $\sum\limits_{k\geq i}m_{k,j} \leq  \sum\limits_{k > i }m_{k,j-1}$ for $i=1,\ldots,n,\overline{n},\ldots,\overline{1}$ and $j=1,\ldots,n-1$, \\
where we set $\sum\limits_{k>i}m_{k,j-1}=0$ for $i=\overline{1}.$
\end{itemize}

\end{Def}
Due to the crystal structure and the equivalence relation above we observe the following remark which helps us to guarantee that operating interchanges with lower decomposition later. Moreover it implies that $\N$ and hence also $\N$ are stable under application of $\e_i$ and $\f_i.$   
\begin{Remark}
Let $m_{i,j}$ be a reduced version of a matrix in  $\mbox{Mat}_{2n \times \mathbb{Z}} (\Z_{\geq 0})$ and  $i\in I.$ 

\begin{itemize}
\item[(i)] If $\f_i$ acts on $m_{i,k}$ then $m_{i,k}>m_{i+1,k-1},$
\item[(i)] if $\e_i$ acts on $m_{i+1,p}$ then $m_{i+1,p}>m_{i,p+1},$

\end{itemize}
where we set $i+1=\overline{i-1}$ if $i\in\{\overline{n},\ldots,\overline{2}\}$ and $n+1=\overline{n}.$ 
\end{Remark}

As mentioned above we also use the lower decomposition rule for the later constructions. Therefore we need an $C$-analog of  Lemma 2.1.

\begin{Lemma}
Let $[M]=(m_{ij})_{{i=1,\ldots,n,\overline{n},\ldots,\overline{1}} \atop {j=0,\ldots, l-1}} \in \mbox{Mat}_{2n \times l} (\Z_{\geq 0})$ be a reduced matrix without zero columns which satisfies condition $(ii)$ of Definition $4.1$.\\ Then we have $$l\leq n.$$   
 \end{Lemma}
\begin{proof}
Assume $l>n$ and consider a  special collection of elements in \\ $\B=\{1,\ldots,n,\bar{n},\ldots\bar{1}\}$: 
\begin{center}
For $k=1,\ldots,n+1$ let $i_k\in \B$ be maximal such that $m_{i_k,k}\neq 0.$
\end{center}
This collection exists because there are no zero columns and $l>n$. Furthermore condition (ii) of Definition 3.1 implies
$$
i_{n+1} \prec i_{n} \prec \ldots \prec i_2 \prec i_1.
$$ 
That means there exists at least one pair $p,q \in \{1,\ldots,n+1\}$ with $p>q$ such that 
$$i_p\in\{1,\ldots,n\} \mbox{ and } i_q=\overline{i_p}.$$ 
Let $p$ be minimal with this property. The minimality of $p$ yields 
$$p-q \leq n-i_p+1.$$    
We assume  $p-q<n-i_p+1,$ namely $p-q=n-i_p+1-j$ for some $j\in \mathbb{N}$.\\
Let us consider the number of elements between $i_p$ and $i_q$:
$$|\{i_p,i_{p-1}, \ldots, i_{q+1},i_{q}\}|=n-i_p-j+2.$$ 
Since
$$(n-i_p-j+2)+(i_p-1)=n-j+1<n+1$$
there is another pair $\hat{p},\hat{q}$ with $\hat{p}>p>\hat{q}$ such that  $i_{\hat{p}}\in\{1,\ldots,n\}$ and $i_{\hat{q}}=\overline{i_{\hat{p}}}$.
If we consider $\hat{p}$ minimal with this property one gets 
$$\hat{p}-\hat{q} \leq n-i_{\hat{p}}+1.$$
If we assume $\hat{p}-\hat{q} < n-i_{\hat{p}}+1$ we can use the same arguments as above.
This way we can inductively conclude that there has to be such a pair with $p-q=n-i_p+1$ which is a contradiction to $m_{i,j}$ being reduced and therefore $l\leq n$.

\qed

\end{proof}

Let us define the $C$-analog of the map $\Phi$ given in section 2. 
$$
\begin{array}{ccccccl}
&\Phi&   :   &\mbox{Mat}_{2n \times \mathbb{Z}} (\Z_{\geq 0})/ \sim&    \rightarrow  &\N / \sim& \\
     &&         &m_{ij}&    \mapsto       &n_{ij},&
     \end{array}              
$$
where we compute $n_{ij}$ as follows: \\

Let $M$ be a reduced version matrix in $\mbox{Mat}_{2n \times \mathbb{Z}} (\Z_{\geq 0})/ \sim$. Then we consider the lower decomposition of $M$:
$$M=M_1+M_2,$$
with $M_1=m_{i,j}^{(1)}$, $M_2=m_{i,j}^{(2)}$ and $M_1\in \N.$
We use exactly the same decomposition as in Section 2.1 with $2n$ rows instead of $n+1$.  Then we move every entry of $M_2$ one column to the left and denote this matrix by $M_2^{(1)}$ and set $$M^{(1)}:=M_1+M_2^{(1)}.$$
Then we decompose $[M^{(1)}]$ and proceed the same way until the iteration becomes stationary and we reach $M^{(k)}=m_{i,j}^{(k)}\in \N.$ Then set
$$n_{i,j}=N:=M^{(k)}.$$
Before we show that this map has the desired properties we state another lemma which will be useful for the proof of the main theorem.  

\begin{Lemma} Let $M=M_1+M_2$ be  the lower decomposition of a matrix in reduced form with $M_1=m_{i,j}^{(1)}$ and $M_2=m_{i,j}^{(2)}$. Then there exists no pair $p,q$ with $p-q=n-i$ such that 
$$m_{i,p}^{(1)}\neq0 \mbox{ and } m_{\overline{i},q}^{(2)}\neq0.$$  
\end{Lemma} 
\begin{proof}
Due to the lower decomposition rule we obtain $M_1 \in \N$ and since $M$ is reduced it`s obvious that $M_1$ is reduced. Now one can use the same arguments as used in the proof of Lemma 4.1.

\qed 

\end{proof}

We define the \emph{compression map} $\kappa$ again as the following composition  
$$
\kappa:=\Psi^{-1} \circ \Phi \circ \Psi 
$$
and show
\begin{Theo} Let $\mathfrak{g}$ be of type $C$. Then, the map
$$
\begin{array}{ccccccl}
&\kappa&   :   &\M&    \rightarrow  &\bigcup\limits_{\lambda\in P,s\in\Z} \mathcal{M}_s(\lambda)& \\
    &&              &M&    \mapsto      &(\Psi^{-1} \circ \Phi \circ \Psi) (M)&
     \end{array}              
$$
defined as above is a morphism of crystals.

\end{Theo}
\begin{proof}
We limit ourselves to prove that sending a reduced matrix $M$ onto $M^{(1)}$ thus defined preserves the crystal structure. This implies inductively that $\Phi$ and hence $\Psi^{-1} \circ \Phi \circ \Psi$ are crystal morphisms. So consider $M=m_{l,j}$ the reduced form of an arbitrary matrix in $\mbox{Mat}_{2n \times \mathbb{Z}} (\Z_{\geq 0})$ and $m_{l,j}=m_{l,j}^{(1)}+m_{l,j}^{(2)}$ its lower decomposition. For an $i \in I$ we have to show:

\begin{itemize}
\item[(i)]$wt(M)=wt(M^{(1)})$,
\item[(ii)]$\ph_i(M)=\ph_i(M^{(1)}),$
\item[(iii)]$\ep_i(M)=\ep_i(M^{(1)})$
\end{itemize}

and that computing $M^{(1)}$ commutes with the Kashiwara operators namely
\begin{center}
\begin{itemize}
\item[(iv)]$(\f_i(M))^{(1)}=\f_i(M^{(1)}),$
\item[(v)]$(\e_i(M))^{(1)}=\e_i(M^{(1)}).$
\end{itemize}
\end{center}

We confine ourselves to prove (ii) and (iv) because the rest follows analogously. Rather we just show (ii) and get (iv) from the $A_n$-case. Set $N=M_1+M_2^{(1)}$ with $N=n_{l,j}$ 
and let $k$ be minimal such that 
$$
\ph_i(M) = \sum\limits_{j \leq k}m_{ij}+m_{\overline{i+1},j-n+i} - \sum\limits_{j <k} m_{i+1,j}+m_{\overline{i},j-n+i}.
$$
In the first step we show $\ph_i(M) \leq \ph_i(N).$ \\
For simplicity we introduce some notation:\\
 For $i\in I,l\in \Z$ and a matrix $M=m_{i,j}$ we set
$$
p_{i,l}(M):=\sum\limits_{j \leq l}m_{i,j}+m_{\overline{i+1},j-n+i} - \sum\limits_{j <l} m_{i+1,j}+m_{\overline{i},j-n+i}.
$$
First of all we look at the case $m_{\overline{i+1},k-n+i}=0$. That yields $m_{i,k}\neq 0$ and furthermore Remark 3.2 says $m_{i,k}>m_{i+1,k-1}$. The fact that $M$ is reduced also implies $m_{\overline{i},k-n+i-1}=0.$  We show that either 
$$
p_{i,k}(N)=p_{i,k}(M)=\ph_i(M)
$$
or
$$
p_{i,k-1}(N)=p_{i,k}(M)=\ph_i(M).
$$ \\
Let us assume that $p_{i,k}(N)<p_{i,k}(M)$. The first case that could yield this is 
$$\sum\limits_{j \leq k}n_{i,j}-\sum\limits_{j < k}n_{i+1,j}<\sum\limits_{j \leq k}m_{i,j}-\sum\limits_{j < k}m_{i+1,j}.$$
But the computation in the proof of Theorem 2.1 shows that in this case we obtain
$$\sum\limits_{j \leq k}m_{i,j}-\sum\limits_{j < k}m_{i+1,j}=\sum\limits_{j \leq k-1}n_{i,j}-\sum\limits_{j < k-1}m_{i+1,j}$$
and since $m_{\overline{i+1},k-n+i}=0,m_{\overline{i},k-n+i-1}=0$ we get
$$\ph_i(M)=p_{i,k-1}(N).$$
The other possibility to get $p_{i,k}(N)<p_{i,k}(M)$ is 
$$\sum\limits_{j \leq k}n_{i,j}-\sum\limits_{j < k}n_{i+1,j}=\sum\limits_{j \leq k}m_{i,j}-\sum\limits_{j < k}m_{i+1,j}$$
and 
$$\sum\limits_{j \leq k}n_{\overline{i+1},j-n+i}-\sum\limits_{j <k}n_{\overline{i},j-n+i}<\sum\limits_{j \leq k}m_{\overline{i+1},j-n+i}-\sum\limits_{j <k}m_{\overline{i},j-n+i}.$$
The first equation implies $0\neq m_{i,k}^{(1)}>m_{i+1,k-1}^{(1)}$ and the inequation yields \\  
$0 \neq m_{\overline{i},k-n+i}^{(2)}>m_{\overline{i+1},k-n+i+1}^{(2)}$.\\
But the existence of $0\neq m_{i,k}^{(1)}$ and $0 \neq m_{\overline{i},k-n+i}^{(2)}$ provides a contradiction to Lemma 4.2.\\
\\
Now we consider the case $m_{\overline{i+1},k-n+i}\neq 0$. Since $M$ is reduced we get $m_{i+1,k}=0$ and in particular $m_{i+1,k}^{(2)}=0$.
Therefore the only chance to have $p_{i,k}(N)<p_{i,k}(M)$ is $m_{i,k}^{(1)}\neq 0$ and $0 \neq m_{\overline{i},k-n+i}^{(2)}>m_{\overline{i+1},k-n+i+1}^{(2)},$ which is again a contradiction to Lemma 4.2. \\ Overall we get 
$$\ph_i(M) \leq \ph_i(N).$$
It remains to show: 
$$\ph_i(M) \geq \ph_i(N).$$
Suppose $\ph_i(M) < \ph_i(N)$, that means there is a $t \in \Z$ with
$$p_{i,t}(N)>\ph_i(M).$$
We distinguish the same cases as above. At first we assume
$$\sum\limits_{j \leq t}n_{i,j}-\sum\limits_{j < t}n_{i+1,j}>\sum\limits_{j \leq t}m_{i,j}-\sum\limits_{j < t}m_{i+1,j}.$$
That is only possible if $m_{i,t+1}^{(2)}> m_{i+1,t}^{(2)}$. In particular we obtain $m_{i,t+1}\neq 0$ and hence $m_{\overline{i},p-n+i}=0.$  
From the $A_n$-case we know that the following inequation holds in this case 
$$\sum\limits_{j \leq t}n_{i,j}-\sum\limits_{j < t}n_{i+1,j} \leq \sum\limits_{j \leq t+1}m_{i,j}-\sum\limits_{j < t+1}m_{i+1,j}.$$  
Combining this with $m_{\overline{i},p-n+i}=0$ and $p_{i,t}(N)>\ph_i(M)$ we see
$$ p_{i,t+1}(M)>\ph_i(M),$$
which is a contradiction to the choice of $k$. \\
Now suppose that 
$$ \sum\limits_{j \leq t}n_{\overline{i+1},j-n+i}-\sum\limits_{j <t}n_{\overline{i},j-n+i}>\sum\limits_{j \leq t}m_{\overline{i+1},j-n+i}-\sum\limits_{j <t}m_{\overline{i},j-n+i}.$$
That means $0 \neq m_{\overline{i+1},t-n+i+1}^{(2)}>m_{\overline{i},t-n+i}^{(2)}$ and
the $A_n$-case implies again  
$$\sum\limits_{j \leq t}n_{\overline{i},j-n+i}-\sum\limits_{j < t}n_{\overline{i},j-n+i} \leq \sum\limits_{j \leq t+1}m_{\overline{i},j-n+i}-\sum\limits_{j < t+1}m_{\overline{i},j-n+i}.$$
If $m_{i,p+1}=0$ we obtain the same contradiction as in the above case. So we assume $m_{i,p+1}\neq0.$ But in order to get $p_{i,t}(N)>\ph_i(M)$ we need $m_{i,p+1}^{(1)}\neq 0$ because otherwise  
$$p_{i,t}(N)\leq p_{i,t+1}(M)\leq\ph_i(M).$$
But $m_{i,p+1}^{(1)}\neq 0$ and $m_{\overline{i+1},t-n+i+1}^{(2)}\neq 0$ provide a contradiction to Lemma 4.2. \\
Hence 
$$\ph_i(M) \geq \ph_i(N).$$
Moreover these arguments also show that either $k$ or $k-1$ is minimal such that either

$$\ph_i(N)=p_{i,k}(N)$$
or
$$\ph_i(N)=p_{i,k-1}(N).$$
Finally Remark 4.2 and the $A_n$-case imply (iv) which finishes our proof.

\qed

\end{proof}

\begin{Example}  For $\mathfrak{g}$ of type $C_3$ we consider the monomial 
$$M=Y_1(0)Y_1(2)Y_1(1)^{-1}Y_1(5)^{-1}Y_1(3)^{-1}Y_1(4)^{-2}Y_2(0)Y_2(3)Y_2(5)^{-2}Y_3(0)Y_3(4).$$
We can write $M$ as 
$$
\begin{array}{ccl}
M &=&X_1(0) X_1(2) X_{\overline{1}}(-2)X_{\overline{1}}(2) X_{\overline{1}}(0)   X_{\overline{1}}(1)^2 X_2(0)X_1(1)X_2(3)X_1(4) X_{\overline{2}}(3)^2 \\
    &&  X_{\overline{1}}(2)^2  X_3(0)X_2(1)X_1(2)X_3(4)X_2(5) X_1(6)
\end{array}
$$
with reduced version matrix
$$[m_{i,j}]=
\begin{pmatrix}           1& 0 &2& 0 & 0   \\ 
                                   1& 1& 0& 1 & 0   \\       
                                   1& 0& 0& 0  & 1  \\
                                   0& 0& 0& 0 & 0   \\
                                   0& 0& 0& 1 & 0   \\  
                                   1& 1& 2& 0 & 0   \\
       \end{pmatrix}.
$$
We observe its lower decomposition by
$$
\begin{pmatrix}           1& 0 &2& 0 & 0   \\ 
                                   1& 1& 0& 1 & 0   \\       
                                   1& 0& 0& 0  & 1  \\
                                   0& 0& 0& 0 & 0   \\
                                   0& 0& 0& 1 & 0   \\  
                                   1& 1& 2& 0 & 0   \\
       \end{pmatrix}
       =
\begin{pmatrix}           1& 0 &1& 0 & 0   \\ 
                                   1& 1& 0& 0 & 0   \\       
                                   1& 0& 0& 0  & 0  \\
                                   0& 0& 0& 0 & 0   \\
                                   0& 0& 0& 0 & 0   \\  
                                   1& 0& 0& 0 & 0   \\
       \end{pmatrix}
       +
       \begin{pmatrix}     0& 0 &1& 0 & 0   \\ 
                                   0& 0& 0& 1 & 0   \\       
                                   0& 0& 0& 0  & 1  \\
                                   0& 0& 0& 0 & 0   \\
                                   0& 0& 0& 1 & 0   \\  
                                   0& 1& 2& 0 & 0   \\
       \end{pmatrix}
$$
and therefore
$$
[m_{i,j}]^{(1)}=
\begin{pmatrix}           1& 1 &1& 0    \\ 
                                   1& 1& 1& 0    \\       
                                   1& 0& 0& 1     \\
                                   0& 0& 0& 0   \\
                                   0& 0& 1& 0   \\  
                                   2& 2& 0& 0   \\
       \end{pmatrix}=[[m_{i,j}]^{(1)}].
$$
One further step yields the desired matrix in $\N$:
$$
\Phi(\Psi(M))=
\begin{pmatrix}           1& 1 &1    \\ 
                                   1& 2& 0    \\       
                                   1& 0& 1      \\
                                   0& 0& 0   \\
                                   0& 1& 0   \\  
                                   4& 0& 0   \\
       \end{pmatrix}
$$
and by application of $\Psi^{-1}$ we get the monomial 
$$
\begin{array}{ccl}
N &=&  X_1(0)X_1(1)X_1(2)X_2(0)X_2(1)^2X_2(2)X_{\overline{2}}(0) X_{\overline{3}}(0)^4 \\
    &=& Y_1(0)Y_2(0)Y_2(1)^2Y_1(3)^{-1}Y_3(1)^{-4}Y_2(1)^4\in  \mathcal{M}_0(3\Lambda_1+2\Lambda_2+2\Lambda_3).
    \end{array}
$$
\end{Example}
Let $B(M)$ be the connected component of $M\in \M$. 
\begin{Cor} For $M \in \M$ and $\Phi(\Psi(M))=:n_{i,j} \in  \N$ we consider $s \in \mathbb{Z}$ maximal such that $n_{i,j}=0$ for all $j<s$ and $i\in \{1,\ldots,n,\overline{n},\ldots\overline{1}\}= \B$. Furthermore for $k=1,\ldots,n$ we define the values:
$$
a_k:=\sum\limits_{i\in \B} n_{i,k+s-1} -\sum\limits_{i\in \B} n_{i,k+s}\geq 0.
$$  
Then we have 
$$ \kappa(M) \in \mathcal{M}_s(\sum\limits_{k=1}^{n}a_k\Lambda_k)$$
and by restricting our morphism to the connected component we get that

$$
\begin{array}{ccccccl}
&\kappa_{|B(M)} &   :   &B(M)&    \rightarrow  & \mathcal{M}_s(\sum\limits_{k=1}^{n}a_k\Lambda_k)& \\
    
     \end{array}              
$$

is a crystal isomorphism.

\end{Cor}

Kim and Shin [7] also gave a realization of the crystal bases in the sense of reversed Young tableaux for Lie algebras of type $C$. In this case they obtained $S(\lambda)$ as the set of all semistandard reversed Young tableaux of shape $\lambda$ with entries $1,\ldots,n,\overline{n},\ldots,\overline{1}$  satisfying some conditions (for details see $[7]$). Moreover Kang, Kim and Shin [3] constructed a morphism between $\mathcal{M}_1(\lambda)$ and those tableaux for $\frak{g}$ of type $C_n$. This is similar to the one in section 2 and can also be generalized to a crystal morphism between arbitrary monomials and tableaux in $S(\lambda)$ via compression. \\ 

Let $M$ be in $\mathcal{M}_1(\lambda)$ for an integral dominant weight $\lambda$ and $m_{i,j}$ its associated reduced matrix. For $i\in \B$ we define again $S(M)$ to be the reversed tableaux with $m_{i,j}$ many $i$`s in the $j$-th row. In order to get a tableaux that satisfies the condition of $S(\lambda)$ we have to apply the rules (al-1) and (al-2) which are due to $[3]$. If we denote by $[S(M)]$ the reversed tableaux we obtain from  $S(M)$ by applying those rules we can state the following
\begin{Prop}$[3]$
The map
$$
\begin{array}{ccccccl}
&\Omega&  :   &M_1(\lambda)&  \rightarrow  &S(\lambda)& \\
&&                         &M&    \mapsto          &[S(M)]&
\end{array}
$$

is a crystal isomorphism.

\end{Prop}
As in the $A_n$-case we continue this morphism to $\mathcal{M}_s(\lambda)$.
Let $M\in \mathcal{M}_s(\lambda)$ be a monomial and $m_{ij}$ the associated reduced matrix in $\N$. We set $S(M)$ to be the semistandard reversed tableaux with $m_{ij}$-many $i$ entries in the $j-s+1$-st row and get the morphism  
$$\begin{array}{ccccccl}
&\Omega&  :   &\bigcup\limits_{\lambda\in P,s\in \mathbb{Z}} \mathcal{M}_s(\lambda)&  \rightarrow  &\bigcup\limits_{\lambda \in P} S(\lambda)& \\
&&                         &M&    \mapsto          &[S(M)].&
\end{array}
$$
If we combine this result  with Theorem 4.1 we get a morphism between Nakajima monomials and tableaux: 
\begin{Cor}
The map
$$
\begin{array}{ccccccl}
&\Omega \circ \kappa &  :   &\M&  \rightarrow  &\bigcup\limits_{\lambda \in P}S(\lambda)& \\

\end{array}
$$

is a crystal morphism.
\end{Cor}

Let us consider an example:
\begin{Example} For $\mathfrak{g}$ of type $C_3$ consider the monomial $$M=Y_2(2)^2Y_2(1)^{-1}Y_3(0)Y_1(0)Y_3(3)^{-1}.$$
Via compression we get $\Phi(\Psi(M))$ by 
$$\begin{pmatrix}         1& 0 & 1& 0  \\ 
                                     1& 1&  0& 0  \\       
                                    1& 0&   0& 0  \\
                                    0& 1&  0& 0 \\
                                    1& 0&  0& 0 \\  
                                    1& 0&  0& 0 \\
       \end{pmatrix}.
  $$
This yields the following tableau:

$$S(\kappa(M))={\young(::::1,:::2{{\overline{3}}},123{{\overline{2}}}{{\overline{1}}})}\stackrel{(\mbox{al-1})}{=}{\young(::::2,:::3{{\overline{3}}},123{{\overline{3}}}{{\overline{2}}})}=[S(\kappa(M))]\in S(\lambda).$$ 

\end{Example}

\section{Insertion scheme for monomials in type $C$}

In this section we define a $C$-analogue of the bumping rule for Nakajima monomials given in section 3. Let $M_1, M_2 \in \mbox{Mat}_{2n \times \mathbb{Z}} (\Z_{\geq 0})$ be two reduced version matrices of monomials in $\M$. As in the $A_n$-case we need to associate a matrix $M_1*M_2 \in \mbox{Mat}_{2n \times \mathbb{Z}} (\Z_{\geq 0})$ to $M_1\otimes M_2$. In order to assure that $M_1*M_2$ is in reduced form we need more zero columns between $M_1$ and $M_2$ in this case. Namely we insert $n$ zero-columns and define 
$$M_1*M_2= 
\begin{pmatrix}                          
 & &  0 & \ldots & 0 & &  \\
 & &  0 & & 0 & &   \\
 &M_2  & \vdots & & \vdots &M_1 & \\
 & &  0 & & 0 & &   \\
 & &  0 & \ldots& 0 & &

\end{pmatrix} \in \mbox{Mat}_{2n \times \mathbb{Z}} (\Z_{\geq 0}).
$$
With the tensor product rule and the same arguments as in Proposition 3.1 we observe:

\begin{Prop}
The map
$$
\begin{array}{ccccl} 
 &\mbox{Mat}_{2n \times \mathbb{Z}} (\Z_{\geq 0})/ \sim \otimes  \mbox{ Mat}_{2n \times \mathbb{Z}} (\Z_{\geq 0})/ \sim& \rightarrow &\mbox{Mat}_{2n \times \mathbb{Z}} (\Z_{\geq 0})/ \sim& \\
 
 &M_1\otimes M_2& \mapsto &M_1*M_2& 
 \end{array}
 $$ is a crystal morphism.
\end{Prop}
Furthermore for $M_1,M_2\in \M$ we can define the bumping $M_2 \rightarrow M_1$ via compression analogously to the $A_n$-case:
$$ M_1 \rightarrow M_2 := \Psi^{-1}(\Phi(\Psi(M_1)*\Psi(M_2))),$$
where $\Psi$ is the crystal isomorphism between $\M$ and $\mbox{Mat}_{2n \times \mathbb{Z}} (\Z_{\geq 0})/ \sim$ and \\ $\Phi: \mbox{ Mat}_{2n \times \mathbb{Z}} (\Z_{\geq 0})/ \sim \rightarrow \N/\sim$ the matrix compression. Theorem 4.1 and Proposition 5.1 imply 
\begin{Theo} Let $\mathfrak{g}$ be of type C. Then, the map
$$
\begin{array}{ccccl} 
&\M \otimes \M & \rightarrow &\bigcup\limits_{\lambda \in P,s\in\Z}\mathcal{M}_s(\lambda)& \\
  &M_1 \otimes M_2 & \mapsto & M_1 \rightarrow M_2&   

\end{array}
$$
is a crystal morphism. 
\end{Theo}
Kim and Shin [7],[8] also defined a bumping rule for reversed tableaux in type $C$. Therefore it is natural to compare $[S(M_1 \rightarrow M_2)]$ and $[S(M_1)] \rightarrow [S(M_2)]$ as in section 3, where $[S(M_1 \rightarrow M_2)], [S(M_1)]$ and $[S(M_2)]$ are the corresponding tableaux in $\bigcup\limits_{\lambda} S(\lambda)$  due to Corollary 4.2. Theorem 5.1 together with Corollary 4.2 imply again
$$  
[S(M_1 \rightarrow M_2)]=[S(M_1)] \rightarrow [S(M_2)].
$$

\newpage
\section*{Acknowledgements}
The author would like to thank Professor Peter Littelmann for the adoption into this subject and Professor Seok-Jin Kang for a short but inspiring discussion the results.
 
\section*{References}

\begin{trivlist}\setlength\labelwidth{1em}
\item[{[1]}] Kang, S.-J., Kim, S.-J., Lee, H., Shin, D.-U.: \textsl{Young wall realization of crystal bases for classical Lie algebras.}
Trans. Amer. Math. Soc. \textbf{356}, 2349-2378  (2004)
\item[{[2]}]  Kang, S.-J., Kim, S.-J. Shin, D.-U.: \textsl{Monomial realization of crystal bases for special linear Lie algebras.} J. Algebra \textbf{274}, 629-642 (2004)
\item[{[3]}] Kang, S.-J., Kim, J.-A., Shin, D.-U.: \textsl{Crystal bases for quantum classical algebras and Nakajima`s monomials.} Publ. Res. Inst. Math. Sci. \textbf{40}, 758-791 (2004)
\item[{[4]}] Kashiwara, M.: \textsl{On crystal bases of the q-analogue of universal enveloping algebras.} Duke Math. J. \textbf{63}, 465-516 (1991)
\item[{[5]}] Kashiwara, M.: \textsl{Realizations of crystals.} Contemp. Math. \textbf{325}, 133-139 (2003)
\item[{[6]}] Kashiwara, M., Nakashima, T.: \textsl{Crystal graphs for representations of the q-analogue of classical
Lie algebras.} J. Algebra \textbf{165}, 295-345 (1994)
\item[{[7]}]  Kim, J.-A., Shin, D.-U. \textsl{Correspondence between Young walls and Young tableaux realizations of crystal bases for the
classical Lie algebras.} J. Algebra \textbf{282}, 728-757 (2004) 
\item[{[8]}]  Kim, J.-A., Shin, D.-U. \textsl{Insertion schemes for the classical Lie algebras.} Communications in Algebra \textbf{32}, 3139-3167  (2004)
\item[{[9]}] Littelmann, P. \textsl{Paths and root operators in representation theory.} Annals of Math. \textbf{142}, 449-525  (1995) 
\item[{[10]}] Littelmann, P. \textsl{Characters of representations and paths in $\mathfrak{h}_\R^*.$ } Proceedings of Symposia in Pure Mathematics \textbf{61}, 29-49  (1997) 
\item[{[11]}] Nakajima H., \textsl{t-analogs of q-characters of quantum affine algebras of type $A_n$ and $D_n$.} Contemp. Math. \textbf{325}, 141-160 (2003)

\end{trivlist}

\end{document}